\documentclass[12pt]{article}%
\usepackage{amsfonts}
\usepackage{amsmath,amssymb}
\usepackage[applemac]{inputenc}
\usepackage{amsmath,amssymb,fullpage}
\usepackage{color}
\usepackage{amsmath}
\usepackage{amssymb}
\usepackage{graphicx}%
\providecommand{\U}[1]{\protect\rule{.1in}{.1in}}
\newtheorem{theorem}{Theorem}

\newtheorem{definition}[theorem]{Definition}
\newtheorem{example}[theorem]{Example}

\newtheorem{lemma}[theorem]{Lemma}

\newtheorem{problem}[theorem]{Problem}

\begin{document}

\title{Model Uncertainty Stochastic Mean-Field Control}
\author{Nacira Agram$^{1,2}$ and Bernt Øksendal$^{1}$}
\date{20 June 2018}
\maketitle

\begin{abstract}
\noindent We consider the problem of optimal control of a mean-field
stochastic differential equation (SDE) under model uncertainty. The model
uncertainty is represented by ambiguity about the law $\mathcal{L}(X(t))$ of
the state $X(t)$ at time $t$. For example, it could be the law $\mathcal{L}%
_{\mathbb{P}}(X(t))$ of $X(t)$ with respect to the given, underlying
probability measure $\mathbb{P}$. This is the classical case when there is no
model uncertainty. But it could also be the law $\mathcal{L}_{\mathbb{Q}%
}(X(t))$ with respect to some other probability measure $\mathbb{Q}$ or, more
generally, any random measure $\mu(t)$ on $\mathbb{R}$ with total mass
$1$.\newline

\noindent We represent this model uncertainty control problem as a
\emph{stochastic differential game} of a mean-field related type SDE with two
players. The control of one of the players, representing the uncertainty of
the law of the state, is a measure-valued stochastic process $\mu(t)$ and the
control of the other player is a classical real-valued stochastic process
$u(t)$. This optimal control problem with respect to random probability
processes $\mu(t)$ in a non-Markovian setting is a new type of stochastic
control problems that has not been studied before. By constructing a new
Hilbert space $\mathcal{M}$ of measures, we obtain a sufficient and a
necessary maximum principles for Nash equilibria for such games in the general
nonzero-sum case, and for saddle points in zero-sum games.\newline

\noindent As an application we find an explicit solution of the problem of
optimal consumption under model uncertainty of a cash flow described by a
mean-field related type SDE.

\end{abstract}

\footnotetext[1]{Department of Mathematics, University of Oslo, P.O. Box 1053
Blindern, N--0316 Oslo, Norway.\newline Email: \texttt{naciraa@math.uio.no,
oksendal@math.uio.no.}
\par
This research was carried out with support of the Norwegian Research Council,
within the research project Challenges in Stochastic Control, Information and
Applications (STOCONINF), project number 250768/F20.}

\footnotetext[2]{University of Biskra, Algeria.}

\paragraph{MSC(2010):}

\noindent60H05, 60H20, 60J75, 93E20, 91G80, 91B70.

\paragraph{Keywords:}

\noindent Mean-field stochastic differential equation; measure-valued optimal
control; model uncertainty; stochastic differential game; stochastic maximum
principle; operator-valued backward stochastic differential equation; optimal
consumption of a mean-field cash flow under model uncertainty.

\section{Introduction}

\noindent There are many ways of introducing model uncertainty. For example,
in recent works of Øksendal and Sulem \cite{OS2}, \cite{OS}, \cite{OS3}, the
underlaying probability measure is not given a priori and there can be a
family of possible probability measures to choose from.

\noindent The aim of this paper is to study stochastic optimal control under
model uncertainty of a mean-field related type SDE driven by Brownian motion
and an independent Poisson random measure. The model uncertainty is
represented by ambiguity about the law $\mathcal{L}(X(t))$ of the state $X(t)$
at time $t$. For example, it could be the law $\mathcal{L}_{\mathbb{P}}(X(t))$
of $X(t)$ with respect to the given, underlying probability measure
$\mathbb{P}$. This is the classical case when there is no model uncertainty.
But it could also be the law $\mathcal{L}_{\mathbb{Q}}(X(t))$ with respect to
some other probability measure $\mathbb{Q}$ or, more generally, any random
measure $\mu(t)$ on $\mathbb{R}$ with total mass $1$.\newline We represent
this model uncertainty control problem as a \emph{stochastic differential
game} of a mean-field related type SDE with two players. The control of one of
the players, representing the uncertainty of the law of the state, is a
measure-valued stochastic process $\mu(t)$, and the control of the other
player is a classical real-valued stochastic process $u(t)$. We penalize
$\mu(t)$ for being far away from the law $\mathcal{L}_{\mathbb{P}}(X(t))$ with
respect to the original probability measure $\mathbb{P}$. This leads to a new
type of mean-field stochastic control problems in which the control is random
measure-valued stochastic process $\mu(t)$ on $\mathbb{R}$. \newline

\noindent To the best of our knowledge this type of problem has not been
studied before. By constructing a new Hilbert space $\mathcal{M}$ of measures,
we obtain sufficient and necessary maximum principles for Nash equilibria for
such games in the general nonzero-sum case, and saddle points for zero-sum
games. As an application we find an explicit solution of the problem of
optimal consumption under model uncertainty of a cash flow described by a
mean-field related type SDE.\newline

\noindent Mean-field games problems were first studied by Lasry and Lions
\cite{LL} and Lions in \cite{Lions} has proved the differentiability of
functions of measures defined on a Wasserstein metric space $\mathcal{P}_{2}$
by using the lifting technics. Since then this type of problems has gained a
lot attention, we can for example refer to Carmona \textit{et al}
\cite{carmona}, \cite{carmona1}, Buckdahn \textit{et al} \cite{BLP},
Bensoussan \textit{et al} \cite{B}, Bayraktar \textit{et al} \cite{BCP}, Corso
and Pham \cite{CP}, Djehiche and Hamadene \cite{DH}, Pham and Wei \cite{PW}
and Agram \cite{A}.\newline

\section{A weighted Sobolev space of random measures}

\noindent In this section, we as in Agram and Øksendal \cite{AO} construct a
Hilbert space $\mathcal{M}$ of random measures on $\mathbb{R}$. It is simpler
to work with than the Wasserstein metric space that has been used by many
authors previously.

\begin{definition}
{(Weighted Sobolev spaces of measures)} \noindent For $k=0, 1, 2, ...$ let
$\mathcal{\tilde{M}}^{(k)}$ denote the set of random measures $\mu$ on
$\mathbb{R}$ such that%
\begin{equation}
\mathbb{E[}%
%TCIMACRO{\tint _{\mathbb{R}}}%
%BeginExpansion
{\textstyle\int_{\mathbb{R}}}
%EndExpansion
|\hat{\mu}(y)|^{2}|y|^{k}e^{-y^{2}}dy]<\infty, \label{eq1}%
\end{equation}
where
\begin{equation}
\hat{\mu}(y)={{ {\textstyle\int_{\mathbb{R}}} }}e^{ixy}d\mu(x) \label{eq2}%
\end{equation}
is the Fourier transform of the measure $\mu$. \noindent If $\mu,\eta
\in\mathcal{\tilde{M}}^{(k)}$ we define the inner product $\left\langle
\mu,\eta\right\rangle _{\mathcal{\tilde{M}}^{(k)}}$ by%
\begin{equation}
\left\langle \mu,\eta\right\rangle _{\mathcal{\tilde{M}}^{(k)}}=\mathbb{E[}
{\textstyle\int_{\mathbb{R}}} \operatorname{Re}(\overline{\hat{\mu}}%
(y)\hat{\eta}(y))|y|^{k}e^{-y^{2}}dy], \label{eq3}%
\end{equation}
where, in general, $\operatorname{Re}(z)$ denotes the real part of the complex
number $z$, and $\bar{z}$ denotes the complex conjugate of $z$. \noindent The
norm $|| \cdot||_{\tilde{\mathcal{M}}^{(k)}}$ associated to this inner product
is given by
\begin{equation}
\label{eq4}\left\Vert \mu\right\Vert _{\mathcal{\tilde{M}}^{(k)}}%
^{2}=\left\langle \mu,\mu\right\rangle _{\mathcal{\tilde{M}}^{(k)}%
}=\mathbb{E[}%
%TCIMACRO{\tint _{\mathbb{R}}}%
%BeginExpansion
{\textstyle\int_{\mathbb{R}}}
%EndExpansion
|\hat{\mu}(y)|^{2}|y|^{k}e^{-y^{2}}dy]\text{.}%
\end{equation}
\noindent The space $\mathcal{\tilde{M}}^{(k)}$ equipped with the inner
product $\left\langle \mu,\eta\right\rangle _{\mathcal{\tilde{M}}^{(k)}}$ is a
pre-Hilbert space. \noindent We let $\mathcal{M}^{(k)}$ denote the completion
of this pre-Hilbert space. \noindent We denote by $\mathcal{M}^{(k)}_{0}$ the
set of all deterministic elements of $\mathcal{M}^{(k)}$. \noindent For $k=0$
we write $\mathcal{M}^{(0)}=\mathcal{M}$ and $\mathcal{M}^{(0)}_{0}%
=\mathcal{M}_{0}$.
\end{definition}

\noindent There are several advantages with working with this Hilbert space
$\mathcal{M}$, compared to the Wasserstein metric space:

\begin{itemize}
\item Our space of measures is easier to work with.

\item A Hilbert space has a useful stronger structure than a metric space.

\item The Wasserstein metric space $\mathcal{P}_{2}$ deals only with
probability measures with finite second moment, while our Hilbert space deals
with any (random) measure satisfying \eqref{eq1}.

\item With this norm we have the following useful estimate:

\begin{lemma}
Let $X^{(1)}$ and $X^{(2)}$ be two random variables in $L^{2}(\mathbb{P})$.
Then
\[%
\begin{array}
[c]{lll}%
\left\Vert \mathcal{L}(X^{(1)})-\mathcal{L}(X^{(2)})\right\Vert _{\mathcal{M}%
_{0}}^{2} & \leq & \sqrt{\pi}\mathbb{E}[(X^{(1)}-X^{(2)})^{2}]\text{.}%
\end{array}
\]

\end{lemma}

We refer to \cite{AO} for a proof.
\end{itemize}

\noindent Let us give some examples of measures:

\begin{example}
[Measures]\-

\begin{enumerate}
\item Suppose that $\mu=\delta_{x_{0}}$, the unit point mass at $x_{0}%
\in\mathbb{R}$. Then $\delta_{x_{0}}\in\mathcal{M}_{0}$ and
\[
{%
%TCIMACRO{\tint _{\mathbb{R}}}%
%BeginExpansion
{\textstyle\int_{\mathbb{R}}}
%EndExpansion
}e^{ixy}d\mu(x)=e^{ix_{0}y},
\]

and hence
\[%
\begin{array}
[c]{lll}%
\left\Vert \mu\right\Vert _{\mathcal{M}_{0}}^{2} & =%
%TCIMACRO{\tint _{\mathbb{R}}}%
%BeginExpansion
{\textstyle\int_{\mathbb{R}}}
%EndExpansion
|e^{ix_{0}y}|^{2}e^{-y^{2}}dy & <\infty\text{.}%
\end{array}
\]

\item Suppose $d\mu(x)=f(x)dx$, where $f\in L^{1}(\mathbb{R})$. Then $\mu
\in\mathcal{M}_{0}$ and by Riemann-Lebesque lemma, $\hat{\mu}(y)\in
C_{0}(\mathbb{R})$, i.e. $\hat{\mu}$ is continuous and $\hat{\mu
}(y)\rightarrow0$ when $|y|\rightarrow\infty$. In particular, $|\hat{\mu}|$ is
bounded on $\mathbb{R}$ and hence%
\[%
\begin{array}
[c]{lll}%
\left\Vert \mu\right\Vert _{\mathcal{M}_{0}}^{2} & =%
%TCIMACRO{\tint _{\mathbb{R}}}%
%BeginExpansion
{\textstyle\int_{\mathbb{R}}}
%EndExpansion
|\hat{\mu}(y)|^{2}e^{-y^{2}}dy & <\infty\text{.}%
\end{array}
\]

\item Suppose that $\mu$ is any finite positive measure on $\mathbb{R}$. Then
$\mu\in\mathcal{M}^{(k)}_{0}$ for all $k$, because
\[%
\begin{array}
[c]{lll}%
|\hat{\mu}(y)| & \leq%
%TCIMACRO{\tint _{\mathbb{R}}}%
%BeginExpansion
{\textstyle\int_{\mathbb{R}}}
%EndExpansion
d\mu(y)=\mu(\mathbb{R}) & <\infty\text{, for all }y\text{,}%
\end{array}
\]
and hence%
\[%
\begin{array}
[c]{lll}%
\left\Vert \mu\right\Vert _{\mathcal{M}^{(k)}_{0}}^{2} & =%
%TCIMACRO{\tint _{\mathbb{R}}}%
%BeginExpansion
{\textstyle\int_{\mathbb{R}}}
%EndExpansion
|\hat{\mu}(y)|^{2}|y|^{k}e^{-y^{2}}dy \leq\mu^{2}(\mathbb{R}) \int
_{\mathbb{R}} |y|^{k} e^{-y^{2}}dy & <\infty\text{.}%
\end{array}
\]

\item Next, suppose $x_{0}=x_{0}(\omega)$ is random. Then $\delta
_{x_{0}(\omega)}$ is a random measure in $\mathcal{M}$. Similarly, if
$f(x)=f(x,\omega)$ is random, then $d\mu(x,\omega)=f(x,\omega)dx$ is a random
measure in $\mathcal{M}$.\newline
\end{enumerate}
\end{example}

\subsection{t-absolute continuity and t-derivative of the law process}

\noindent Let $(%
%TCIMACRO{\U{3a9} }%
%BeginExpansion
\Omega
%EndExpansion
,\mathcal{F},\mathbb{P})$ be a given probability space with filtration
$\mathbb{F}=(\mathcal{F}_{t})_{t\geq0}$ generated by a one-dimensional
Brownian motion $B$ and an independent Poisson random measure $N(dt,d\zeta)$.
Let $\nu(d\zeta)dt$ denote the Lévy measure of $N$, and let $\tilde
{N}(dt,d\zeta)$ denote the compensated Poisson random measure $N(dt,d\zeta
)-\nu(d\zeta)dt$.

\noindent Suppose that $X(t)=X_{t}$ is an Itô-Lévy process of the form
\begin{equation}%
\begin{cases}
dX_{t}=\alpha(t)dt+\beta(t)dB(t)+%
%TCIMACRO{\tint _{\mathbb{R}_{0}}}%
%BeginExpansion
{\textstyle\int_{\mathbb{R}_{0}}}
%EndExpansion
\gamma(t,\zeta)\tilde{N}(dt,d\zeta);\quad t\in\lbrack0,T],\\
X_{0}=x\in\mathbb{R},
\end{cases}
\label{eq2.1}%
\end{equation}
where $\alpha,\beta$ and $\gamma$ are bounded predictable processes.

\noindent Let $\varphi\in C^{2}$. Then under appropriate conditions on the
coefficients, we get by the Itô\ formula
\begin{equation}
\mathbb{E}[\varphi(X_{t+h})]-\mathbb{E}[\varphi(X_{t})]=\mathbb{E}[{%
%TCIMACRO{\tint _{t}^{t+h}}%
%BeginExpansion
{\textstyle\int_{t}^{t+h}}
%EndExpansion
}A\varphi(X_{s})ds], \label{eq2.3}%
\end{equation}
where
\[%
\begin{array}
[c]{lll}%
A\varphi(X_{s}) & = & \alpha(s)\varphi^{\prime}(X_{s})+\tfrac{1}{2}\beta
^{2}(s)\varphi^{\prime\prime}(X_{s})\\
&  & +%
%TCIMACRO{\tint _{\mathbb{R}_{0}}}%
%BeginExpansion
{\textstyle\int_{\mathbb{R}_{0}}}
%EndExpansion
\{\varphi(X_{s}+\gamma(s,\zeta))-\varphi(X_{s})-\varphi^{\prime}(X_{s}%
)\gamma(s,\zeta)\}\nu(d\zeta).
\end{array}
\]
In particular, if
\[
\varphi(x)=\varphi_{y}(x):=\exp(ixy);\quad y\in\mathbb{R},
\]
then
\[%
\begin{array}
[c]{ll}%
A\varphi_{y}(X_{s}) & =(iy\alpha(s)-\tfrac{1}{2}\beta^{2}(s)y^{2}\\
& \text{ \ \ }+%
%TCIMACRO{\tint _{\mathbb{R}_{0}}}%
%BeginExpansion
{\textstyle\int_{\mathbb{R}_{0}}}
%EndExpansion
\left\{  \exp(i\gamma(s,\zeta)y)-1-iy\gamma(s,\zeta)\right\}  \nu
(d\zeta))\varphi_{y}(X_{s}),
\end{array}
\]
for all $y\in\mathbb{R}$.

\begin{definition}
[Law process]From now on we use the notation
\[
M_{t}:=M(t):=\mathcal{L}(X_{t});\quad0\leq t\leq T
\]
for the law process $\mathcal{L}(X_{t})$ of $X_{t}=X(t)$ with respect to
$\mathbb{P}$.
\end{definition}

\begin{lemma}
(i) The map $t\mapsto M_{t}:[0,T]\rightarrow\mathcal{M}_{0}$ is absolutely
continuous, and the derivative%

\[
M^{\prime}(t):=\frac{d}{dt}M(t)
\]
exists for all $t$.\newline

(ii) There exists a constant $C<\infty$ such that
\begin{equation}
||M^{\prime}(t)||_{\mathcal{M}_{0}} \leq C ||M(t)||_{\mathcal{M}^{(4)}_{0}}
\text{ for all } t \in[0,T]; M(t) \in\mathcal{M}^{(4)}_{0}.
\end{equation}

\end{lemma}

\noindent{Proof.} (i) \quad Let $0\leq t<t+h\leq T$. Then by \eqref{eq2} and
\eqref{eq4} we get
\begin{align}
&  \left\Vert M_{t+h}-M_{t}\right\Vert _{\mathcal{M}_{0}}^{2}=%
%TCIMACRO{\tint _{\mathbb{R}}}%
%BeginExpansion
{\textstyle\int_{\mathbb{R}}}
%EndExpansion
|\hat{M}_{t+h}(y)-\hat{M}_{t}(y)|^{2}e^{-y^{2}}dy\nonumber\\
&  =%
%TCIMACRO{\tint _{\mathbb{R}}}%
%BeginExpansion
{\textstyle\int_{\mathbb{R}}}
%EndExpansion
|%
%TCIMACRO{\tint _{\mathbb{R}}}%
%BeginExpansion
{\textstyle\int_{\mathbb{R}}}
%EndExpansion
e^{ixy}d\mathcal{L}(X_{t+h})-%
%TCIMACRO{\tint _{\mathbb{R}}}%
%BeginExpansion
{\textstyle\int_{\mathbb{R}}}
%EndExpansion
e^{ixy}d\mathcal{L}(X_{t})(x)|^{2}e^{-y^{2}}dy\nonumber\\
&  =%
%TCIMACRO{\tint _{\mathbb{R}}}%
%BeginExpansion
{\textstyle\int_{\mathbb{R}}}
%EndExpansion
|\mathbb{E}[\varphi_{y}(X_{t+h})]-\mathbb{E}[\varphi_{y}(X_{t})]|^{2}%
e^{-y^{2}}dy. \label{2.8}%
\end{align}
The last equality holds by using that for any bounded function $\psi$ we have%
\[%
\begin{array}
[c]{lll}%
\mathbb{E}[\psi(X)] & = &
%TCIMACRO{\tint _{\mathbb{R}}}%
%BeginExpansion
{\textstyle\int_{\mathbb{R}}}
%EndExpansion
\psi(x)d\mathcal{L}(X)(x).
\end{array}
\]
By $\left(  \ref{eq2.3}\right)  $, we obtain%

\begin{align}
\left\Vert M_{t+h}-M_{t}\right\Vert _{\mathcal{M}_{0}}^{2}  &  ={%
%TCIMACRO{\tint _{\mathbb{R}}}%
%BeginExpansion
{\textstyle\int_{\mathbb{R}}}
%EndExpansion
}|\mathbb{E}[{%
%TCIMACRO{\tint _{t}^{t+h}}%
%BeginExpansion
{\textstyle\int_{t}^{t+h}}
%EndExpansion
}A\varphi_{y}(X(s))ds]|^{2}e^{-y^{2}}dy\nonumber\\
&  \leq{%
%TCIMACRO{\tint _{\mathbb{R}}}%
%BeginExpansion
{\textstyle\int_{\mathbb{R}}}
%EndExpansion
}({%
%TCIMACRO{\tint _{t}^{t+h}}%
%BeginExpansion
{\textstyle\int_{t}^{t+h}}
%EndExpansion
}\mathbb{E}[\left\vert A\varphi_{y}(X_{s})\right\vert ]ds)^{2}e^{-y^{2}}dy\leq
C_{1}\text{ }h^{2}, \label{eq9}%
\end{align}
for some constant $C_{1}$ which does not depend on $t$ and $h$.\newline We
have proved that for different $t$ and $t+h$, $\left\Vert M_{t+h}%
-M_{t}\right\Vert _{\mathcal{M}_{0}}^{2}\leq C$ $h^{2}$ and it is easy to see
that this holds for every finite disjoint partition of the interval $[0,T]$.
Thus we get that $t \mapsto M(t)$ is absolutely continuous, and the derivative
$M^{\prime}(t)=\frac{d}{dt}M(t)$ exists for all $t$.\newline

(ii) This follows from \eqref{eq9}, using that the coefficients $\alpha
,\beta,\gamma$ are bounded and that
\begin{equation}
\mathbb{E}[|A_{\varphi_{y}}(X_{s})|]\leq const.y^{2}|\mathbb{E}[\exp
(iyX_{s})]|\leq const.y^{2}|\widehat{M}_{s}(y)|.
\end{equation}
.\hfill$\square$ \newline
%\bigskip
From the lemma above we conclude the following:

\begin{lemma}
If $X_{t}$ is an Itô-Lévy process as in \eqref{eq2.1}, then the derivative
$M_{s}^{\prime}:=\frac{d}{ds}M_{s}$ exists in $\mathcal{M}_{0}$ for a.a. $s$,
and we have
\[
M_{t}=M_{0}+%
%TCIMACRO{\tint _{0}^{t}}%
%BeginExpansion
{\textstyle\int_{0}^{t}}
%EndExpansion
M_{s}^{\prime}ds;\quad t\geq0.
\]

\end{lemma}

\noindent In the following we will apply this to the solutions $X(t)$ of the
mean-field related type SDEs we consider below.

\begin{example}
\-

\begin{description}
\item[(a)] Suppose that $X(t)=B(t)$ with $B(0)=0$. Then
\[
d\mathcal{L}(X(t))(x)=\tfrac{1}{\sqrt{2\pi t}}\exp(-\tfrac{x^{2}}{2t})dx,
\]
i.e. $\mathcal{L}(X(t))$\ has a density $\tfrac{1}{\sqrt{2\pi t}}\exp
(-\tfrac{x^{2}}{2t}).$ Therefore $\tfrac{d}{dt}\mathcal{L}(X(t))$ is a measure
with density%
\[
\tfrac{d}{dt}\tfrac{1}{\sqrt{2\pi t}}\exp(-\tfrac{x^{2}}{2t})=(\tfrac{x^{2}%
-t}{2t^{2}})(\tfrac{1}{\sqrt{2\pi t}}\exp(-\tfrac{x^{2}}{2t})).
\]

\item[(b)] Suppose $X(t)=N(t)$, a Poisson process with intensity $\bar
{\lambda}$. Then for $k=1,2,...$ we have
\[
\mathbb{P}(N(t)=k)=\tfrac{e^{-\bar{\lambda}t}(\bar{\lambda}t)^{k}}{k!}%
\]
and hence
\[
\tfrac{d}{dt}\mathbb{P}(N(t)=k)=\tfrac{1}{k!}(\bar{\lambda}e^{-\bar{\lambda}%
t}(\lambda t)^{k-1}\{k-\bar{\lambda}t\}).
\]

\end{description}
\end{example}

\section{Preliminaries}

\noindent We will recall some concepts and spaces which will be used on the sequel.

\noindent The probability $\mathbb{P}$ is a reference probability measure. We
introduce two smaller filtrations $\mathbb{G}^{(i)}\mathbb{=}(\mathcal{G}%
_{t}^{(i)})_{t\geq0}$ such that $\mathcal{G}_{t}^{(i)}\subseteq\mathcal{F}%
_{t}$, for $i=1,2$ and for all $t\geq0.$ These filtrations represent the
information available to player number $i$ at time $t$. \newline

\subsection{Some basic concepts from Banach space theory}

\noindent Since we deal with measures defined on an Hilbert space
$\mathcal{M}$, we need the Fréchet derivative to differentiate functions of
measures. Let $\mathcal{X},\mathcal{Y}$ be two Banach spaces with norms
$\Vert\cdot\Vert_{\mathcal{X}},\Vert\cdot\Vert_{\mathcal{Y}}$, respectively,
and let $F:\mathcal{X}\rightarrow\mathcal{Y}$.

\begin{itemize}
\item We say that $F$ has a directional derivative (or Gâteaux derivative) at
$v\in\mathcal{X}$ in the direction $w\in\mathcal{X}$ if
\[
D_{w}F(v):=\lim_{\varepsilon\rightarrow0}\frac{1}{\varepsilon}(F(v+\varepsilon
w)-F(v))
\]
exists in $\mathcal{Y}$.

\item We say that $F$ is Fréchet differentiable at $v\in\mathcal{X}$ if there
exists a continuous linear map $A:\mathcal{X}\rightarrow\mathcal{Y}$ such
that
\[
\lim_{\substack{h\rightarrow0\\h\in\mathcal{X}}}\frac{1}{\Vert h\Vert
_{\mathcal{X}}}\Vert F(v+h)-F(v)-A(h)\Vert_{\mathcal{Y}}=0.
\]
In this case we call $A$ the \textit{gradient} (or Fréchet derivative) of $F$
at $v$ and we write
\[
A=\nabla_{v}F.
\]

\item If $F$ is Fréchet differentiable at $v$ with Fréchet derivative
$\nabla_{v}F$, then $F$ has a directional derivative in all directions
$w\in\mathcal{X}$ and
\[
D_{w}F(v):=\left\langle \nabla_{v}F,w\right\rangle =\nabla_{v}F(w)=\nabla
_{v}Fw.
\]

\end{itemize}

\noindent In particular, note that if $F$ is a linear operator, then
$\nabla_{v}F=F$ for all $v$.\newline

\subsection{Spaces}

\noindent Throughout this work, we will use the following spaces:

%\item[b)] Throughout this work, we will use the following spaces:
%\end{description}

\begin{itemize}
\item $\mathcal{S}^{2}$ is the set of ${\mathbb{R}}$-valued $\mathbb{F}%
$-adapted càdlàg processes $(X(t))_{t\in\lbrack0,T]}$ such that
\[
{\Vert X\Vert}_{\mathcal{S}^{2}}^{2}:={\mathbb{E}}[\sup_{t\in\lbrack
0,T]}|X(t)|^{2}]~<~\infty\;,
\]

\item $\mathbb{L}^{2}$ is the set of ${\mathbb{R}}$-valued $\mathbb{F}%
$-predictable processes $(Q(t))_{t\in\lbrack0,T]}$ such that
\[
\Vert Q\Vert_{\mathbb{L}^{2}}^{2}:={\mathbb{E}}[%
%TCIMACRO{\tint _{0}^{T}}%
%BeginExpansion
{\textstyle\int_{0}^{T}}
%EndExpansion
|Q(t)|^{2}dt]<~\infty\;.
\]

\item $L^{2}(\mathcal{F}_{t})$ is the set of ${\mathbb{R}}$-valued square
integrable $\mathcal{F}_{t}$-measurable random variables.

\item $\mathbb{L}_{\nu}^{2}$ is the set of $\mathbb{F}$-predictable processes
$R:[0,T]\times\mathbb{R}_{0}\times\Omega\rightarrow\mathbb{R}$ such that
\[
||R||_{\mathbb{L}_{\nu}^{2}}^{2}:={\mathbb{E}}[%
%TCIMACRO{\tint _{\mathbb{R}_{0}}}%
%BeginExpansion
{\textstyle\int_{\mathbb{R}_{0}}}
%EndExpansion
|R(t,\zeta)|^{2}\nu(d\zeta)dt]~<~\infty\;.
\]

\item In general, for any given filtration $\mathbb{H}$, we say that the
measure-valued process $\mu(t)=\mu(t,\omega):[0,T]\times\Omega\rightarrow
\mathcal{M}$ is adapted to $\mathbb{H}$ if $\mu(t)(V)$ is $\mathbb{H}$-adapted
for all Borel sets $V\subseteq%
%TCIMACRO{\U{211d} }%
%BeginExpansion
\mathbb{R}
%EndExpansion
$. Let $\mathbb{M}_{\mathbb{G}}=\mathbb{M}_{\mathbb{G}^{1}}$ be a given set of
$\mathcal{M}$-valued, $\mathbb{G}^{1}=(\mathcal{G}_{t}^{1})_{t\geq0}%
$-predictable, stochastic processes $\mu(t)$. We call $\mathbb{M}_{\mathbb{G}%
}$ the set of admissible measure-valued control processes $\mu(\cdot)$.

\item $\mathbb{M}_{0}$ is the set of $t$-differentiable $\mathcal{M}_{0}%
$-valued processes $m(t); t \in[0,T]$.\newline If $m \in\mathbb{M}_{0}$ we put
$m^{\prime}(t)=\frac{d}{dt}m(t)$.

\item Let $\mathcal{A}_{\mathbb{G}}=\mathcal{A}_{\mathbb{G}^{2}}$ be a given
set of real-valued, $\mathbb{G}^{2}=(\mathcal{G}_{t}^{2})_{t\geq0}%
$-predictable, stochastic processes $u(t)$ required to have values in a given
convex subset $\mathcal{U}$ of $%
%TCIMACRO{\U{211d} }%
%BeginExpansion
\mathbb{R}
%EndExpansion
$. We call $\mathcal{A}_{\mathbb{G}}$ the set of admissible real-valued
control processes $u(\cdot)$.

\item $\mathcal{R}$ is the set of measurable functions $r:\mathbb{R}%
_{0}\rightarrow\mathbb{R}.$

\item $C_{a}([0,T],\mathcal{M}_{0})$ denotes the set of absolutely continuous
functions $m:[0,T]\rightarrow\mathcal{M}_{0}.$

\item $\mathbb{K}$ is the set of bounded linear functionals $K:\mathcal{M}%
_{0}\rightarrow\mathbb{R}$ equipped with the operator norm
\[
||K||_{\mathbb{K}}:=\sup_{m\in\mathcal{M}_{0},||m||_{\mathcal{M}_{0}}\leq
1}|K(m)|.
\]

\item $\mathcal{S}_{\mathbb{K}}^{2}$ is the set of $\mathbb{F}$-adapted càdlàg
processes $p:[0,T]\times\Omega\mapsto\mathbb{K}$ such that
\[
||p||_{\mathcal{S}_{\mathbb{K}}}^{2}:=\mathbb{E}[\sup_{t\in\lbrack
0,T]}||p(t)||_{\mathbb{K}}^{2}]<\infty.
\]

\item $\mathbb{L}_{\mathbb{K}}^{2}$ is the set of $\mathbb{F}$-predictable
processes $q:[0,T]\times\Omega\mapsto\mathbb{K}$ such that
\[
||q||_{\mathbb{L}_{\mathbb{K}}^{2}}^{2}:=\mathbb{E}[%
%TCIMACRO{\tint _{0}^{T}}%
%BeginExpansion
{\textstyle\int_{0}^{T}}
%EndExpansion
||q(t)||_{\mathbb{K}}^{2}dt]<\infty.
\]

\item $\mathbb{L}_{\nu,\mathbb{K}}^{2}$ is the set of $\mathbb{F}$-predictable
processes $r:[0,T]\times\mathbb{R}_{0}\times\Omega\mapsto\mathbb{K}$ such
that
\[
||r||_{\mathbb{L}_{\nu,\mathbb{K}}^{2}}^{2}:=\mathbb{E}[%
%TCIMACRO{\tint _{0}^{T}}%
%BeginExpansion
{\textstyle\int_{0}^{T}}
%EndExpansion%
%TCIMACRO{\tint _{\mathbb{R}_{0}}}%
%BeginExpansion
{\textstyle\int_{\mathbb{R}_{0}}}
%EndExpansion
||r(t,\zeta)||_{\mathbb{K}}^{2}\nu(d\zeta)dt]<\infty.
\]

\end{itemize}

\section{The model uncertainty stochastic optimal control problem}

\noindent As pointed out in the Introduction, there are several ways to
represent model uncertainty in a stochastic system. In this paper, we are
interested in systems governed by controlled mean-field related type SDE
$X^{\mu,u}(t)=X(t)\in\mathcal{S}^{2}$ on the form%

\begin{equation}
\left\{
\begin{array}
[c]{lll}%
dX(t) & = & b\left(  t,X(t),\mu(t),u(t)\right)  dt+\sigma\left(
t,X(t),\mu(t),u(t)\right)  dB(t)\\
&  & +%
%TCIMACRO{\tint _{\mathbb{R}_{0}}}%
%BeginExpansion
{\textstyle\int_{\mathbb{R}_{0}}}
%EndExpansion
\gamma\left(  t,X(t),\mu(t),u(t),\zeta\right)  \tilde{N}(dt,d\zeta);\text{
}t\in\left[  0,T\right]  ,\\
X\left(  0\right)  & = & x\in%
%TCIMACRO{\U{211d} }%
%BeginExpansion
\mathbb{R}
%EndExpansion
.
\end{array}
\right.  \label{sde}%
\end{equation}
The functions
\[%
\begin{array}
[c]{llll}%
b(t,x,\mu,u) & =b(t,x,\mu,u,\omega) & :\left[  0,T\right]  \times%
%TCIMACRO{\U{211d} }%
%BeginExpansion
\mathbb{R}
%EndExpansion
\times\mathcal{M}\times\mathcal{U}\times\Omega & \rightarrow%
%TCIMACRO{\U{211d} }%
%BeginExpansion
\mathbb{R}
%EndExpansion
,\\
\sigma(t,x,\mu,u) & =\sigma(t,x,\mu,u,\omega) & :\left[  0,T\right]  \times%
%TCIMACRO{\U{211d} }%
%BeginExpansion
\mathbb{R}
%EndExpansion
\times\mathcal{M}\times\mathcal{U}\times\Omega & \rightarrow%
%TCIMACRO{\U{211d} }%
%BeginExpansion
\mathbb{R}
%EndExpansion
,\\
\gamma(t,x,\mu,u,\zeta) & =\gamma(t,x,\mu,u,\zeta,\omega) & :\left[
0,T\right]  \times%
%TCIMACRO{\U{211d} }%
%BeginExpansion
\mathbb{R}
%EndExpansion
\times\mathcal{M}\times\mathcal{U}\times%
%TCIMACRO{\U{211d} }%
%BeginExpansion
\mathbb{R}
%EndExpansion
_{0}\times\Omega & \rightarrow%
%TCIMACRO{\U{211d} }%
%BeginExpansion
\mathbb{R}
%EndExpansion
,
\end{array}
\]
are supposed to be Lipschitz on $x\in%
%TCIMACRO{\U{211d} }%
%BeginExpansion
\mathbb{R}
%EndExpansion
$, uniformly with respect to $t$ and $\omega$ for given $u\in\mathcal{U}$ and
$\mu\in\mathcal{M}$. Then by e.g. Theorem $1.19$ in Øksendal and Sulem
\cite{osjump}, we have existence and uniqueness of the solution of $X(t)$. We
may regard \eqref{sde} as a perturbed version of the mean-field equation
\begin{equation}
\left\{
\begin{array}
[c]{lll}%
dX(t) & = & b\left(  t,X(t),\mathcal{L}(X(t)),u(t)\right)  dt+\sigma\left(
t,X(t),\mathcal{L}(X(t)),u(t)\right)  dB(t)\\
&  & +%
%TCIMACRO{\tint _{\mathbb{R}_{0}}}%
%BeginExpansion
{\textstyle\int_{\mathbb{R}_{0}}}
%EndExpansion
\gamma\left(  t,X(t),\mathcal{L}(X(t)),u(t),\zeta\right)  \tilde{N}%
(dt,d\zeta);\text{ }t\in\left[  0,T\right]  ,\\
X\left(  0\right)  & = & x\in%
%TCIMACRO{\U{211d} }%
%BeginExpansion
\mathbb{R}
%EndExpansion
.
\end{array}
\right.  \label{sde2}%
\end{equation}
For example, we could have $\mu(t)=\mathcal{L}_{\mathbb{Q}}(X(t))$ for some
probability measure $\mathbb{Q}\neq\mathbb{P}$. \newline Thus the model
uncertainty is represented by an uncertainty about what law $\mu(t)$ is
influencing the coefficients of the system, and we are penalising the laws
that are far away from $\mathcal{L}(X(t))$. See the application in Section
5.\newline

\noindent Let us consider a performance functional of the form
\begin{equation}%
\begin{array}
[c]{lll}%
J\left(  \mu,u\right)  & = & \mathbb{E[}g\left(  X(T),M(T)\right)  +%
%TCIMACRO{\tint _{0}^{T}}%
%BeginExpansion
{\textstyle\int_{0}^{T}}
%EndExpansion
\ell\left(  s,X\left(  s\right)  ,M(s),\mu(s),u\left(  s\right)  \right)  ds],
\end{array}
\label{performance}%
\end{equation}
where $\ell(t,x,m,\mu,u)=\ell(t,x,m,\mu,u,\omega):\left[  0,T\right]  \times%
%TCIMACRO{\U{211d} }%
%BeginExpansion
\mathbb{R}
%EndExpansion
\times\mathcal{M}_{0}\times\mathcal{M}\times\mathcal{U}\times\Omega\rightarrow%
%TCIMACRO{\U{211d} }%
%BeginExpansion
\mathbb{R}
%EndExpansion
$ and $g:%
%TCIMACRO{\U{211d} }%
%BeginExpansion
\mathbb{R}
%EndExpansion
\times\mathcal{M}_{0}\times\Omega\rightarrow%
%TCIMACRO{\U{211d} }%
%BeginExpansion
\mathbb{R}
%EndExpansion
$ are given functions.

\noindent For fixed $x,m,\mu,u$ we assume that $\ell\left(  s,\cdot\right)  $
is $\mathcal{F}_{s}$-measurable for all $s\in\lbrack0,T]$ and $g(\cdot,\cdot)$
is $\mathcal{F}_{T}$-measurable. We also assume the following integrability
condition
\[
\mathbb{E[}\left\vert g\left(  X(T),M(T)\right)  \right\vert ^{2}+%
%TCIMACRO{\tint _{0}^{T}}%
%BeginExpansion
{\textstyle\int_{0}^{T}}
%EndExpansion
\left\vert \ell\left(  s,X\left(  s\right)  ,M(s),\mu(s),u\left(  s\right)
\right)  \right\vert ^{2}ds]<\infty,
\]
for all $\mu\in\mathbb{M_{\mathbb{G}}}$ and $u\in\mathcal{A}_{\mathbb{G}}.$

\noindent Note that the system $\left(  \ref{sde}\right)  $ and the
performance $\left(  \ref{performance}\right)  $ are not Markovian. However,
recently a dynamic programming approaches to mean-field stochastic control
problems have been introduced. See e.g. Bayraktar \textit{et al} \cite{BCP}
and Pham and Wei \cite{PW}. In this paper we will use an approach based on a
suitably modified stochastic maximum principle, which also works in partial
information settings. \newline

\noindent In the next section we study a stochastic differential game of two
players, where one of the players is solving an optimal measure-valued control
problem of the type described above, while the other player is solving a
classical real-valued stochastic control problem. To the best of our knowledge
this type of stochastic differential game has not been studied before.\newline

\subsection{Nonzero-sum games}

\noindent We now proceed to a nonzero-sum maximum principle.\newline We
consider the $\mathbb{R} \times\mathcal{M}_{0}$-valued process $(X(t),M(t))$
where $M(t)=\mathcal{L}(X(t))$, where $X(t)$ is given by \eqref{sde} and
\begin{equation}
dM(t)=\beta(M(t))dt; \quad M(0) \in\mathcal{M}_{0} \text{ given },
\end{equation}
where $\beta$ is the operator on $\mathbb{M}_{0}$ given by
\begin{equation}
\beta(m(t))=m^{\prime}(t).
\end{equation}

\noindent The cost functionals are assumed to be on the form%
\begin{equation}%
\begin{array}
[c]{lll}%
J_{i}\left(  \mu,u\right)  & = & \mathbb{E}[g_{i}\left(  X(T),M(T)\right) \\
&  & +%
%TCIMACRO{\tint _{0}^{T}}%
%BeginExpansion
{\textstyle\int_{0}^{T}}
%EndExpansion
\ell_{i}\left(  s,X\left(  s\right)  ,M(s),\mu(s),u\left(  s\right)  \right)
ds];\text{ for }i=1,2,
\end{array}
\label{perf}%
\end{equation}
where $M(s):=\mathcal{L}(X(s))$ and the functions
\[%
\begin{array}
[c]{llll}%
\ell_{i}(t,x,m,\mu,u) & =\ell_{i}(t,x,m,\mu,u,\omega) & :\left[  0,T\right]
\times%
%TCIMACRO{\U{211d} }%
%BeginExpansion
\mathbb{R}
%EndExpansion
\times\mathcal{M}_{0}\times\mathcal{M}\times\mathcal{U}\times\Omega &
\rightarrow%
%TCIMACRO{\U{211d} }%
%BeginExpansion
\mathbb{R}
%EndExpansion
,\\
g_{i}(x,m) & =g_{i}(x,m,\omega) & :%
%TCIMACRO{\U{211d} }%
%BeginExpansion
\mathbb{R}
%EndExpansion
\times\mathcal{M}_{0}\times\Omega & \rightarrow%
%TCIMACRO{\U{211d} }%
%BeginExpansion
\mathbb{R}
%EndExpansion
,
\end{array}
\]
are continuously differentiable with respect to $x,u$ and admit Fréchet
derivatives with respect to $m$ and $\mu$.

\begin{problem}
We consider the general nonzero-sum stochastic game to find $\left(  \mu
^{\ast},u^{\ast}\right)  \in\mathbb{M}_{\mathbb{G}}\mathcal{\times
A}_{\mathbb{G}}$ such that%
\[%
\begin{array}
[c]{ll}%
J_{1}(\mu,u^{\ast})\leq J_{1}(\mu^{\ast},u^{\ast}), & \text{for all }\mu
\in\mathbb{M}_{\mathbb{G}},\\
J_{2}(\mu^{\ast},u)\leq J_{2}(\mu^{\ast},u^{\ast}), & \text{for all }%
u\in\mathcal{A}_{\mathbb{G}}.
\end{array}
\]
The pair $\left(  \mu^{\ast},u^{\ast}\right)  $ is called a \textit{Nash
equilibrium}.
\end{problem}

\begin{definition}
\textbf{{(The Hamiltonian)}} For $i=1,2$ we define the Hamiltonian
\[
H_{i}:[0,T]\times\mathbb{R}\times\mathcal{M}_{0} \times\mathcal{M}%
\times\mathcal{U}\times\mathbb{R}\times\mathbb{R}\times\mathcal{R}\times
C_{a}([0,T],\mathcal{M}_{0})\rightarrow\mathbb{R}%
\]
by%
\begin{equation}%
\begin{array}
[c]{ll}%
H_{i}(t,x,m,\mu,u,p_{i}^{0},q_{i}^{0},r_{i}^{0}(\cdot),p_{i}^{1}) & =\ell
_{i}(t,x,m,\mu,u)+p_{i}^{0}b(t,x,\mu,u)+q_{i}^{0}\sigma(t,x,\mu,u)\\
& +%
%TCIMACRO{\tint _{\mathbb{R}_{0}}}%
%BeginExpansion
{\textstyle\int_{\mathbb{R}_{0}}}
%EndExpansion
r_{i}^{0}(\zeta)\gamma\left(  t,x,\mu,u,\zeta\right)  \nu(d\zeta)+\left\langle
p_{i}^{1},\beta(m)\right\rangle .
\end{array}
\label{haml}%
\end{equation}

\end{definition}

\noindent We assume that $H_{i}$ is continuously differentiable with respect
to $x,u$ and admits Fréchet derivatives with respect to $m$ and $\mu$.\newline
For $u\in\mathcal{A}_{\mathbb{G}},\mu\in\mathbb{M}_{\mathbb{G}}$ with
corresponding solution $X=X^{\mu,u}$, define $p_{i}^{0}=p_{i}^{0,\mu,u}%
,q_{i}^{0}=q_{i}^{0,\mu,u}$ and $r_{i}^{0}=r_{i}^{0,\mu,u\text{ }}$ and
$p_{i}^{1}=p_{i}^{1,\mu,u},q_{i}^{1}=q_{i}^{1,\mu,u}$ and $r_{i}^{1}%
=r_{i}^{1,\mu,u\text{ }}$for $i=1,2$ by the following set of adjoint equations:

\begin{itemize}
\item The real-valued BSDE in the unknown $(p_{i}^{0},q_{i}^{0},r_{i}^{0}%
)\in\mathcal{S}^{2}\times\mathbb{L}^{2}\times\mathbb{L}_{\nu}^{2}$ is given by%

\begin{equation}
\left\{
\begin{array}
[c]{lll}%
dp_{i}^{0}(t) & = & -\frac{\partial H_{i}}{\partial x}(t)dt+q_{i}%
^{0}(t)dB(t)+\int_{%
%TCIMACRO{\U{211d} }%
%BeginExpansion
\mathbb{R}
%EndExpansion
_{0}}r_{i}^{0}(t,\zeta)\tilde{N}(dt,d\zeta);\quad t\in\lbrack0,T],\\
p_{i}^{0}(T) & = & \frac{\partial g_{i}}{\partial x}(X(T),M(T)),
\end{array}
\right.  \label{eqp0}%
\end{equation}

\item and the operator-valued BSDE in the unknown $(p_{i}^{1},q_{i}^{1}%
,r_{i}^{1})\in\mathcal{S}_{\mathbb{K}}^{2}\times\mathbb{L}_{\mathbb{K}}%
^{2}\times\mathbb{L}_{\nu,\mathbb{K}}^{2}$ is given by
\begin{equation}
\left\{
\begin{array}
[c]{lll}%
dp_{i}^{1}(t) & = & - \nabla_{m}H_{i}(t) dt+q_{i}^{1}(t)dB(t)+\int_{
\mathbb{R} _{0}}r_{i}^{1}(t,\zeta)\tilde{N}(dt,d\zeta);\quad t\in
\lbrack0,T],\\
p^{1}(T) & = & \nabla_{m}g_{i}(X(T),M(T)),
\end{array}
\right.  \label{eqp1}%
\end{equation}

\end{itemize}

\noindent where $H_{i}(t)=H_{i}(t,X\left(  t\right)  ,M(t),\mu(t),u\left(
t\right)  ,p_{i}^{0}(t),q_{i}^{0}(t),r_{i}^{0}(t,\cdot),p_{i}^{1}(t))$
etc.\newline We remark that the BSDEs $\left(  \ref{eqp0}\right)  $
%and $\left(  \ref{eqp1}\right)  $ are
is linear, so whenever knowing the Hamiltonian $H_{i}$ and the function
$g_{i}$, we can get a solution explicitly. To remind the reader of this
solution formula, let us consider the solution $(P,Q,R)\in\mathcal{S}%
^{2}\times\mathbb{L}^{2}\times\mathbb{L}_{\nu}^{2}$ of the linear BSDE
\begin{equation}
\left\{
\begin{array}
[c]{ll}%
dP(t) & =-[\varphi(t)+\alpha(t)P(t)+\beta(t)Q(t)+\int_{%
%TCIMACRO{\U{211d} }%
%BeginExpansion
\mathbb{R}
%EndExpansion
_{0}}\phi(t,\zeta)R(t,\zeta)\nu(d\zeta)]dt\\
& \text{ \ \ \ \ \ \ \ \ \ }+Q(t)dB(t)+\int_{%
%TCIMACRO{\U{211d} }%
%BeginExpansion
\mathbb{R}
%EndExpansion
_{0}}R(t,\zeta)\tilde{N}(dt,d\zeta);\text{ }t\in\left[  0,T\right]  \text{,}\\
P(T) & =\theta\in L^{2}(\mathcal{F}_{T}).
\end{array}
\right.  \label{lobsde}%
\end{equation}
Here $\varphi,\alpha,\beta$ and $\phi$ are bounded predictable processes with
$\phi$ is assumed to be an $%
%TCIMACRO{\U{211d} }%
%BeginExpansion
\mathbb{R}
%EndExpansion
$-valued process defined on $[0,T]\times%
%TCIMACRO{\U{211d} }%
%BeginExpansion
\mathbb{R}
%EndExpansion
_{0}\times\Omega$. Then it is well-known (see e.g. Theorem $1.7$ in Øksendal
and Sulem \cite{OS3}) that the component $P(t)$ of the solution of equation
$\left(  \ref{lobsde}\right)  $ can be written in closed form as follows:
\begin{equation}%
\begin{array}
[c]{lll}%
P(t) & =\mathbb{E[}\theta\tfrac{\Gamma(T)}{\Gamma(t)}+%
%TCIMACRO{\tint _{t}^{T}}%
%BeginExpansion
{\textstyle\int_{t}^{T}}
%EndExpansion
\tfrac{\Gamma(s)}{\Gamma(t)}\varphi(s)|\mathcal{F}_{t}];\text{ } & t\in\left[
0,T\right]  \text{,}%
\end{array}
\label{solutio}%
\end{equation}
where $\Gamma(t)\in\mathcal{S}^{2}$ is the solution of the linear SDE with
jumps%
\begin{equation}
\left\{
\begin{array}
[c]{lll}%
d\Gamma(t) & = & \Gamma(t^{-})[\alpha(t)dt+\beta(t)dB(t)+\int_{%
%TCIMACRO{\U{211d} }%
%BeginExpansion
\mathbb{R}
%EndExpansion
_{0}}\phi(t,\zeta)\tilde{N}(dt,d\zeta)];\text{ }t\in\left[  0,T\right]
\text{,}\\
\Gamma(0) & = & 1.
\end{array}
\right.  \label{gama}%
\end{equation}
For notational convenience, we will employ the following short hand notations%
\[%
\begin{array}
[c]{lll}%
\hat{H}_{1}(t) & = & H_{1}(t,\hat{X}(t),\hat{M}(t),\hat{\mu}(t),\hat
{u}(t),\hat{p}_{1}^{0}(t),\hat{q}_{1}^{0}(t),\hat{r}_{1}^{0}(t,\cdot),\hat
{p}_{1}^{1}(t)),\\
\check{H}_{1}(t) & = & H_{1}(t,\hat{X}(t),\hat{M}(t),\mu(t),\hat{u}(t),\hat
{p}_{1}^{0}(t),\hat{q}_{1}^{0}(t),\hat{r}_{1}^{0}(t,\cdot),\hat{p}_{1}%
^{1}(t)),\\
\bar{H}_{2}(t) & = & H_{2}(t,\hat{X}(t),\hat{M}(t),\hat{\mu}(t),\hat
{u}(t),\hat{p}_{2}^{0}(t),\hat{q}_{2}^{0}(t),\hat{r}_{2}^{0}(t,\cdot),\hat
{p}_{2}^{1}(t)),\\
\breve{H}_{2}(t) & = & H_{2}(t,\hat{X}(t),\hat{M}(t),\hat{\mu}(t),u(t),\hat
{p}_{2}^{0}(t),\hat{q}_{2}^{0}(t),\hat{r}_{2}^{0}(t,\cdot),\hat{p}_{2}%
^{1}(t)).
\end{array}
\]
Similar notation is used for the derivatives of $H,\ell,g,b,\sigma,\gamma$
etc. \newline We now state a sufficient theorem for the nonzero-sum games.

\begin{theorem}
[Sufficient nonzero-sum maximum principle]\label{Thm} Let $\left(  \hat{\mu
},\hat{u}\right)  \in\mathbb{M}_{\mathbb{G}}\mathcal{\times A}_{\mathbb{G}}$
with corresponding solutions $\hat{X}$, $(p_{i}^{0},q_{i}^{0},r_{i}^{0})$ and
$(p_{i}^{1},q_{i}^{1},r_{i}^{1})$ of the forward and backward stochastic
differential equations $\left(  \ref{sde}\right)  $, $\left(  \ref{eqp0}%
\right)  $ and $\left(  \ref{eqp1}\right)  $ respectively. Suppose that

\begin{enumerate}
\item (Concavity) The functions
\[%
\begin{array}
[c]{ll}%
(x,m,\mu) & \mapsto H_{1}(t),\\
(x,m,u) & \mapsto H_{2}(t),\\
(x,m) & \mapsto g_{i}(x,m)\text{, for }i=1,2\text{,}%
\end{array}
\]
are concave $\mathbb{P}$.a.s for each $t\in\left[  0,T\right]  $.

\item (Maximum conditions)%
\begin{equation}%
\begin{array}
[c]{lll}%
\mathbb{E}[\hat{H}_{1}(t)|\mathcal{G}_{t}^{(1)}] & = & \underset{\mu
\in\mathbb{M}_{\mathbb{G}}}{ess\text{ }\sup}\mathbb{E}[\check{H}%
_{1}(t)|\mathcal{G}_{t}^{(1)}],
\end{array}
\label{maxQ}%
\end{equation}

\end{enumerate}

and%
\[%
\begin{array}
[c]{lll}%
\mathbb{E}[\bar{H}_{2}(t)|\mathcal{G}_{t}^{(2)}] & = & \underset
{u\in\mathcal{A}_{\mathbb{G}}}{ess\text{ }\sup}\mathbb{E}[\breve{H}%
_{2}(t)|\mathcal{G}_{t}^{(2)}],
\end{array}
\]

$\mathbb{P}$.a.s for each $t\in\left[  0,T\right]  .$

Then $\left(  \hat{\mu},\hat{u}\right)  $ is a Nash equilibrium for our problem.
\end{theorem}

\noindent{Proof.} \quad Let us first prove that $J_{1}(\mu,\hat{u})\leq
J_{1}(\hat{\mu},\hat{u}).$ \newline By the definition of the cost functional
$\left(  \ref{perf}\right)  $ we have for fixed $\hat{u}\in\mathcal{A}%
_{\mathbb{G}}$ and arbitrary $\mu\in\mathbb{M}_{\mathbb{G}}$
\begin{equation}%
\begin{array}
[c]{lll}%
J_{1}(\mu,\hat{u})-J_{1}(\hat{\mu},\hat{u}) & = & I_{1}+I_{2},
\end{array}
\label{j}%
\end{equation}
where
\begin{align*}
&
\begin{array}
[c]{lll}%
I_{1} & = & \mathbb{E}[\int_{0}^{T}\{\check{\ell}_{1}(t)-\hat{\ell}%
_{1}(t)\}dt],
\end{array}
\\
&
\begin{array}
[c]{lll}%
I_{2} & = & \mathbb{E}[\check{g}_{1}(X(T),M(T))-\hat{g}_{1}(\hat{X}(T),\hat
{M}(T))].
\end{array}
\end{align*}
By the definition of the Hamiltonian $\left(  \ref{haml}\right)  $ we have%
\begin{equation}%
\begin{array}
[c]{ll}%
I_{1} & =\mathbb{E[}%
%TCIMACRO{\tint _{0}^{T}}%
%BeginExpansion
{\textstyle\int_{0}^{T}}
%EndExpansion
\check{H}_{1}(t)-\hat{H}_{1}(t)-\hat{p}_{1}^{0}(t)\tilde{b}(t)-\hat{q}_{1}%
^{0}(t)\tilde{\sigma}(t)-{%
%TCIMACRO{\tint _{\mathbb{R}_{0}}}%
%BeginExpansion
{\textstyle\int_{\mathbb{R}_{0}}}
%EndExpansion
}\hat{r}_{1}^{0}(t,\zeta)\tilde{\gamma}(t,\zeta)\nu(d\zeta)-\langle\hat{p}%
_{1}^{1}(t),\tilde{M}^{\prime}(t)\rangle dt],
\end{array}
\label{i1}%
\end{equation}
where $\tilde{b}(t)=\check{b}(t)-\hat{b}(t)$ etc. By the concavity of $g_{1}$
and the terminal values of the BSDEs $\left(  \ref{eqp0}\right)  $, $\left(
\ref{eqp1}\right)  $, we have%
\[%
\begin{array}
[c]{lll}%
I_{2} & \leq\mathbb{E}[\tfrac{\partial g_{1}}{\partial x}(T)\tilde
{X}(T)+\langle\nabla_{m}g_{1}(T),\tilde{M}(T)\rangle] & =\mathbb{E}[\hat
{p}_{1}^{0}(T)\tilde{X}(T)+\langle\hat{p}_{1}^{1}(T),\tilde{M}(T)\rangle].
\end{array}
\]
Applying the Itô formula to $\hat{p}_{1}^{0}(t)\tilde{X}(t)$ and $\langle
\hat{p}_{1}^{1}(t),\tilde{M}(t)\rangle$, we get
\begin{align}
I_{2}  &  \leq\mathbb{E}[\hat{p}_{1}^{0}(T)\tilde{X}(T)+\langle\hat{p}_{1}%
^{1}(T),\tilde{M}(T)\rangle]\nonumber\\
&  =\mathbb{E}[%
%TCIMACRO{\tint _{0}^{T}}%
%BeginExpansion
{\textstyle\int_{0}^{T}}
%EndExpansion
\hat{p}_{1}^{0}(t)d\tilde{X}(t)+%
%TCIMACRO{\tint _{0}^{T}}%
%BeginExpansion
{\textstyle\int_{0}^{T}}
%EndExpansion
\tilde{X}(t)d\hat{p}_{1}^{0}(t)+%
%TCIMACRO{\tint _{0}^{T}}%
%BeginExpansion
{\textstyle\int_{0}^{T}}
%EndExpansion
\hat{q}_{1}^{0}(t)\tilde{\sigma}(t)dt+%
%TCIMACRO{\tint _{0}^{T}}%
%BeginExpansion
{\textstyle\int_{0}^{T}}
%EndExpansion%
%TCIMACRO{\tint _{\mathbb{R}_{0}}}%
%BeginExpansion
{\textstyle\int_{\mathbb{R}_{0}}}
%EndExpansion
\hat{r}_{1}^{0}(t,\zeta)\tilde{\gamma}(t,\zeta)\nu(d\zeta)dt]\nonumber\\
&  +\mathbb{E}[%
%TCIMACRO{\tint _{0}^{T}}%
%BeginExpansion
{\textstyle\int_{0}^{T}}
%EndExpansion
\langle\hat{p}_{1}^{1}(t),d\tilde{M}(t)\rangle+%
%TCIMACRO{\tint _{0}^{T}}%
%BeginExpansion
{\textstyle\int_{0}^{T}}
%EndExpansion
\tilde{M}(t)d\hat{p}_{1}^{1}(t)]\nonumber\\
&  =\mathbb{E}[%
%TCIMACRO{\tint _{0}^{T}}%
%BeginExpansion
{\textstyle\int_{0}^{T}}
%EndExpansion
\hat{p}_{1}^{0}(t)\tilde{b}(t)dt-%
%TCIMACRO{\tint _{0}^{T}}%
%BeginExpansion
{\textstyle\int_{0}^{T}}
%EndExpansion
\tfrac{\partial\hat{H}_{1}}{\partial x}(t)\tilde{X}(t)dt+%
%TCIMACRO{\tint _{0}^{T}}%
%BeginExpansion
{\textstyle\int_{0}^{T}}
%EndExpansion
\hat{q}_{1}^{0}(t)\tilde{\sigma}(t)dt\nonumber\\
&  +%
%TCIMACRO{\tint _{0}^{T}}%
%BeginExpansion
{\textstyle\int_{0}^{T}}
%EndExpansion%
%TCIMACRO{\tint _{\mathbb{R}_{0}}}%
%BeginExpansion
{\textstyle\int_{\mathbb{R}_{0}}}
%EndExpansion
\hat{r}_{1}^{0}(t,\zeta)\tilde{\gamma}(t,\zeta)\nu(d\zeta)dt+%
%TCIMACRO{\tint _{0}^{T}}%
%BeginExpansion
{\textstyle\int_{0}^{T}}
%EndExpansion
\langle\hat{p}_{1}^{1}(t),\tilde{M}^{\prime}\rangle dt\nonumber\\
&  -{\textstyle\int_{0}^{T}}\langle\nabla_{m}\hat{H}_{1}(t),\tilde
{M}(t)\rangle dt], \label{I2}%
\end{align}
where we have used that the $dB(t)$ and $\tilde{N}(dt,d\zeta)$ integrals with
the necessary integrability property are martingales and then have mean zero.
Substituting $\left(  \ref{i1}\right)  $ and $\left(  \ref{I2}\right)  $ in
$\left(  \ref{j}\right)  $, yields%
\begin{align*}
&  J_{1}(\mu,\hat{u})-J_{1}(\hat{\mu},\hat{u})\\
&  \leq\mathbb{E}[%
%TCIMACRO{\tint _{0}^{T}}%
%BeginExpansion
{\textstyle\int_{0}^{T}}
%EndExpansion
\{\check{H}_{1}(t)-\hat{H}_{1}(t)-\tfrac{\partial\hat{H}_{1}}{\partial
x}(t)\tilde{X}(t)-\langle\nabla_{m}\hat{H}_{1}(t),\tilde{M}(t)\rangle\}dt].
\end{align*}
By the concavity of $H_{1}$ and the fact that the process $\mu$ is
$\mathcal{G}_{t}^{(1)}$-adapted, we obtain%
\begin{align}
J_{1}(\mu,\hat{u})-J_{1}(\hat{\mu},\hat{u})  &  \leq\mathbb{E}[%
%TCIMACRO{\tint _{0}^{T}}%
%BeginExpansion
{\textstyle\int_{0}^{T}}
%EndExpansion
\tfrac{\partial\hat{H}_{1}}{\partial\mu}(t)\left(  \mu(t)-\hat{\mu}(t)\right)
dt]\nonumber\\
&  =\mathbb{E}[%
%TCIMACRO{\tint _{0}^{T}}%
%BeginExpansion
{\textstyle\int_{0}^{T}}
%EndExpansion
\mathbb{E}(\tfrac{\partial\hat{H}_{1}}{\partial\mu}(t)\left(  \mu(t)-\hat{\mu
}(t)\right)  |\mathcal{G}_{t}^{(1)})dt]\nonumber\\
&  =\mathbb{E}[%
%TCIMACRO{\tint _{0}^{T}}%
%BeginExpansion
{\textstyle\int_{0}^{T}}
%EndExpansion
\mathbb{E}(\tfrac{\partial\hat{H}_{1}}{\partial\mu}(t)|\mathcal{G}_{t}%
^{(1)})\left(  \mu(t)-\hat{\mu}(t)\right)  dt]\nonumber\\
&  \leq0,\nonumber
\end{align}
where $\frac{\partial\hat{H}_{1}}{\partial\mu}=\nabla_{\mu}\hat{H}_{1}.$ The
last equality holds because of the maximum condition of $\hat{H}_{1}$ at
$\mu=\hat{\mu}$.\newline Similar considerations apply to prove that
$J_{2}(\hat{\mu},u)\leq J_{2}(\hat{\mu},\hat{u})$. For the sake of
completeness, we give details in the Appendix.\hfill$\square$

\noindent We now state and prove a necessary version of the maximum principle.
We assume the following:

\begin{itemize}
\item Whenever $\mu\in\mathbb{M}_{\mathbb{G}}$ $(u\in\mathcal{A}_{\mathbb{G}%
})$ and $\eta\in\mathbb{M}_{\mathbb{G}}$ $(\pi\in\mathcal{A}_{\mathbb{G}})$
are bounded, there exists $\epsilon>0$ such that
\[
\mu+\lambda\eta\in\mathbb{M}_{\mathbb{G}}\text{ }(u+\lambda\pi\in
\mathcal{A}_{\mathbb{G}})\text{, for each }\lambda\in\left[  -\epsilon
,\epsilon\right]  .
\]

\item For each $t_{0}\in\left[  0,T\right]  $ and each bounded $\mathcal{G}%
_{t_{0}}^{(1)}$-measurable random measure $\alpha_{1}$ and $\mathcal{G}%
_{t_{0}}^{(2)}$-measurable random variable $\alpha_{2}$, the process
\begin{equation}
\eta\left(  t\right)  =\alpha_{1}\mathbf{1}_{\left[  t_{0},T\right]
}(t)\text{ } \label{eta}%
\end{equation}
belongs to $\mathbb{M}_{\mathbb{G}}$ and the process%
\[
\pi\left(  t\right)  =\alpha_{2}\mathbf{1}_{\left[  t_{0},T\right]  }(t)
\]
belongs to $\mathcal{A}_{\mathbb{G}}$.\medskip\newline

\begin{definition}
In general, if $K^{u}(t)$ is a process depending on $u$, we define the
differential operator D on $K$ by
\[
DK^{u}(t):=D^{\pi}K^{u}(t)=\tfrac{d}{d\lambda}K^{u+\lambda\pi}(t)|_{\lambda=0}%
\]
whenever the derivative exists.\newline
\end{definition}
\end{itemize}

\noindent The \emph{derivative} of the state $X(t)$ defined by $\left(
\ref{sde}\right)  $ is
\[
DX^{\mu}(t):=\tfrac{d}{d\lambda}X^{\mu+\lambda\eta}|_{\lambda=0}=Z(t)
\]
exists, and is given by
\begin{equation}
\left\{
\begin{array}
[c]{lll}%
dZ\left(  t\right)  & = & [\tfrac{\partial b}{\partial x}\left(  t\right)
Z\left(  t\right)  +\tfrac{\partial b}{\partial\mu}\left(  t\right)
\eta\left(  t\right)  ]dt+[\frac{\partial\sigma}{\partial x}\left(  t\right)
Z\left(  t\right)  +\frac{\partial\sigma}{\partial\mu}\left(  t\right)
\eta\left(  t\right)  ]dB(t)\\
&  & +%
%TCIMACRO{\tint _{\mathbb{R}_{0}}}%
%BeginExpansion
{\textstyle\int_{\mathbb{R}_{0}}}
%EndExpansion
[\frac{\partial\gamma}{\partial x}\left(  t,\zeta\right)  Z(t)+\frac
{\partial\gamma}{\partial\mu}\left(  t,\zeta\right)  \eta\left(  t\right)
]\tilde{N}(dt,d\zeta);\text{ \ \ \ }t\in\left[  0,T\right]  ,\\
Z\left(  0\right)  & = & 0.
\end{array}
\right.  \label{dervz}%
\end{equation}

\noindent We remark that this derivative process is a linear SDE, then by
assuming that $b$, $\sigma$ and $\gamma$ admit bounded partial derivatives
with respect to $x$ and $\mu$, there is a unique solution $Z(t)\in
\mathcal{S}^{2}$ of $\left(  \ref{dervz}\right)  $.\newline We want to prove
that $Z\left(  t\right)  $ is exactly the derivative in $\mathbb{L}%
^{2}(\mathbb{P})$ of $X^{\mu+\lambda\eta}(t)$ with respect to $\lambda$ at
$\lambda=0.$ More precisely, we want to prove the following.

\begin{lemma}%
\begin{equation}
\mathbb{E[}%
%TCIMACRO{\tint _{0}^{T}}%
%BeginExpansion
{\textstyle\int_{0}^{T}}
%EndExpansion
(\tfrac{X^{\mu+\lambda\eta}(t)-X^{\mu}(t)}{\lambda}-Z\left(  t\right)
)^{2}dt]\rightarrow0\text{ as }\lambda\rightarrow0. \label{vraider}%
\end{equation}

\end{lemma}

\noindent{Proof.} For notational convenience, we have here used the simplified
notations
\begin{equation}
\mu^{\lambda}:=\mu+\lambda\eta\label{mul}%
\end{equation}
and by $X^{\mu^{\lambda}}$ we mean the corresponding solution%

\[
X^{\mu^{\lambda}}(t)=x+%
%TCIMACRO{\tint _{0}^{t}}%
%BeginExpansion
{\textstyle\int_{0}^{t}}
%EndExpansion%
%TCIMACRO{\tint _{\mathbb{R}_{0}}}%
%BeginExpansion
{\textstyle\int_{\mathbb{R}_{0}}}
%EndExpansion
\gamma(s,X^{\mu^{\lambda}}(s),\mu^{\lambda}(s),\zeta)\tilde{N}(ds,d\zeta
);\ \ \ t\in\left[  0,T\right]  ,
\]
when assuming that $b=\sigma=0$, and because $u$ is fixed we can omit it.
Then, by the Itô-Lévy isometry, we get
\[%
\begin{array}
[c]{l}%
\mathbb{E}[%
%TCIMACRO{\tint _{0}^{T}}%
%BeginExpansion
{\textstyle\int_{0}^{T}}
%EndExpansion
(\tfrac{X^{\mu^{\lambda}}(t)-X(t)}{\lambda}-Z\left(  t\right)  )^{2}dt]\\
=\mathbb{E[}%
%TCIMACRO{\tint _{0}^{T}}%
%BeginExpansion
{\textstyle\int_{0}^{T}}
%EndExpansion%
%TCIMACRO{\tint _{\mathbb{R}_{0}}}%
%BeginExpansion
{\textstyle\int_{\mathbb{R}_{0}}}
%EndExpansion
\{\tfrac{\gamma(s,X^{\mu^{\lambda}}(s),\mu^{\lambda}(s),\zeta)-\gamma\left(
s,X(s),\mu(s),\zeta\right)  }{\lambda}-\tfrac{\partial\gamma}{\partial
x}\left(  s,\zeta\right)  Z(t)-\frac{\partial\gamma}{\partial\mu}\left(
s,\zeta\right)  \eta\left(  s\right)  \}\tilde{N}(ds,d\zeta))^{2}dt]\\
=\mathbb{E[}\int_{0}^{T}\int_{%
%TCIMACRO{\U{211d} }%
%BeginExpansion
\mathbb{R}
%EndExpansion
_{0}}{\int_{0}^{t}}(\tfrac{\gamma(s,X^{\mu^{\lambda}}(s),\mu^{\lambda
}(s),\zeta)-\gamma\left(  s,X(s),\mu(s),\zeta\right)  }{\lambda}%
-\tfrac{\partial\gamma}{\partial x}\left(  s,\zeta\right)  Z(s)-\frac
{\partial\gamma}{\partial\mu}\left(  s,\zeta\right)  \eta\left(  s\right)
)^{2}\nu(d\zeta)dsdt].
\end{array}
\]
This goes to $0$ when $\lambda$ goes to $0$, by the bounded convergence
theorem and our assumption on $\gamma$.

\hfill$\square$ \newline

\begin{theorem}
\textbf{{(Necessary nonzero-sum maximum principle)}} Let $\left(  \hat{\mu
},\hat{u}\right)  \in\mathbb{M}_{\mathbb{G}}\mathcal{\times A}_{\mathbb{G}}$
with corresponding solutions $\hat{X}$, $(p_{i}^{0},q_{i}^{0},r_{i}^{0})$ and
$(p_{i}^{1},q_{i}^{1},r_{i}^{1})$ of the forward and backward stochastic
differential equations $\left(  \ref{sde}\right)  $ and
$\eqref{eqp0}-\eqref{eqp1}$, with the corresponding derivative process
$\hat{Z}$ given by $\left(  \ref{dervz}\right)  $. Then the following (i) and
(ii) are equivalent:\medskip

\begin{description}
\item[(i)] For all $\mu$, $\eta\in\mathbb{M}_{\mathbb{G}}$ and for all $u$,
$\pi\in\mathcal{A}_{\mathbb{G}}$
\[
\tfrac{d}{d\lambda}J_{1}(\mu+\lambda\eta,u)|_{\lambda=0}=\tfrac{d}{ds}%
J_{2}(\mu,u+s\pi)|_{s=0}=0,
\]

\item[(ii)]
\[
\mathbb{E}[\tfrac{\partial H_{1}}{\partial\mu}(t)|\mathcal{G}_{t}%
^{(1)}]=\mathbb{E}[\tfrac{\partial H_{2}}{\partial u}(t)|\mathcal{G}_{t}%
^{(2)}]=0.\newline%
\]

\end{description}
\end{theorem}

\noindent{Proof.} \quad First note that, by using the linearity of
$\langle\cdot,\cdot\rangle$ and the fact that the Fréchet derivative of a
linear operator is the same operator, we get, by interchanging the order of
the derivatives $\frac{d}{dt}$ and $\nabla_{m}$, that%

\[
\nabla_{m}\langle p_{1}^{1}(t),\frac{d}{dt}m\rangle=\langle p_{1}%
^{1}(t),\nabla_{m}\frac{d}{dt}m\rangle=\langle p_{1}^{1}(t),\frac{d}{dt}%
\nabla_{m}(m)\rangle=\langle p_{1}^{1}(t),\frac{d}{dt}(\cdot)\rangle,
\]
and hence
\[%
\begin{array}
[c]{lll}%
\langle\nabla_{m}\langle p_{1}^{1}(t),\frac{d}{dt}m\rangle,DM(t)\rangle &
=\langle p_{1}^{1}(t),\frac{d}{dt}DM(t)\rangle & =\left\langle p_{1}%
^{1}(t),DM^{\prime}(t)\right\rangle
\end{array}
\]
Also, note that
%$\left\langle p_{1}^{1}(t),DM^{\prime}(t)\right\rangle=p_{1}^{1}(t)DM^{\prime}(t)$ where%
\[%
\begin{array}
[c]{ll}%
dDM(t) & =DM^{\prime}(t)dt.
\end{array}
\]
Assume that (i) holds. Using the definition of $J_{1}\left(  \ref{perf}%
\right)  $, we get
\begin{align*}
0  &  =\tfrac{d}{d\lambda}J_{1}(\mu+\lambda\eta,u)|_{\lambda=0}\\
&  =\mathbb{E[}%
%TCIMACRO{\tint _{0}^{T}}%
%BeginExpansion
{\textstyle\int_{0}^{T}}
%EndExpansion
\{\tfrac{\partial\ell_{1}}{\partial x}\left(  t\right)  Z\left(  t\right)
+\langle\nabla_{m}\ell_{1}(t),DM\left(  t\right)  \rangle+\tfrac{\partial
\ell_{1}}{\partial\mu}\left(  t\right)  \eta\left(  t\right)  \}dt\\
&  +\tfrac{\partial g_{1}}{\partial x}\left(  T\right)  Z\left(  T\right)
+\langle\nabla_{m}g_{1}\left(  T\right)  ,DM(T)\rangle].
\end{align*}
Hence, by the definition \eqref{haml} of $H_{1}$, we have
\begin{align}
0  &  =\tfrac{d}{d\lambda}J_{1}(\mu+\lambda\eta,u)|_{\lambda=0}\nonumber\\
&  =\mathbb{E[}%
%TCIMACRO{\tint _{0}^{T}}%
%BeginExpansion
{\textstyle\int_{0}^{T}}
%EndExpansion
\{\tfrac{\partial H_{1}}{\partial x}(t)-p_{1}^{0}(t)\tfrac{\partial
b}{\partial x}(t)-q_{1}^{0}(t)\tfrac{\partial\sigma}{\partial x}(t)-%
%TCIMACRO{\tint _{\mathbb{R}_{0}}}%
%BeginExpansion
{\textstyle\int_{\mathbb{R}_{0}}}
%EndExpansion
r_{1}^{0}(t,\zeta)\tfrac{\partial\gamma}{\partial x}\left(  t,\zeta\right)
\nu(d\zeta)\}Z(t)dt\nonumber\\
&  +{\textstyle\int_{0}^{T}}\langle\nabla_{m}H_{1}\left(  t\right)  ,DM\left(
t\right)  \rangle dt\nonumber\\
&  -{\textstyle\int_{0}^{T}}\langle p_{1}^{1}(t),DM^{\prime}\left(  t\right)
\rangle dt+%
%TCIMACRO{\tint _{0}^{T}}%
%BeginExpansion
{\textstyle\int_{0}^{T}}
%EndExpansion
\{\tfrac{\partial H_{1}}{\partial\mu}(t)-p_{1}^{0}(t)\tfrac{\partial
b}{\partial\mu}(t)\nonumber\\
&  -q_{1}^{0}(t)\tfrac{\partial\sigma}{\partial\mu}(t)-%
%TCIMACRO{\tint _{\mathbb{R}_{0}}}%
%BeginExpansion
{\textstyle\int_{\mathbb{R}_{0}}}
%EndExpansion
r_{1}^{0}(t,\zeta)\tfrac{\partial\gamma}{\partial\mu}\left(  t,\zeta\right)
\nu(d\zeta)\}\eta(t)dt+p_{1}^{0}(T)Z(T)+\langle p_{1}^{1}(T),DM(T)\rangle].
\label{dj}%
\end{align}
Applying now the Itô formula to both $p_{1}^{0}Z$ and $\langle p_{1}%
^{1},DM\rangle$, we get
\begin{align}
&  \mathbb{E}[p_{1}^{0}(T)Z(T)+\langle p_{1}^{1}(T),DM(T)\rangle]\nonumber\\
&  =\mathbb{E}[%
%TCIMACRO{\tint _{0}^{T}}%
%BeginExpansion
{\textstyle\int_{0}^{T}}
%EndExpansion
p_{1}^{0}(t)dZ(t)+%
%TCIMACRO{\tint _{0}^{T}}%
%BeginExpansion
{\textstyle\int_{0}^{T}}
%EndExpansion
Z(t)dp_{1}^{0}(t)+%
%TCIMACRO{\tint _{0}^{T}}%
%BeginExpansion
{\textstyle\int_{0}^{T}}
%EndExpansion
q_{1}^{0}(t)(\tfrac{\partial\sigma}{\partial x}\left(  t\right)  Z\left(
t\right)  +\tfrac{\partial\sigma}{\partial\mu}\left(  t\right)  \eta\left(
t\right)  )dt\nonumber\\
&  +%
%TCIMACRO{\tint _{0}^{T}}%
%BeginExpansion
{\textstyle\int_{0}^{T}}
%EndExpansion%
%TCIMACRO{\tint _{\mathbb{R}_{0}}}%
%BeginExpansion
{\textstyle\int_{\mathbb{R}_{0}}}
%EndExpansion
r_{1}^{0}(t,\zeta)(\tfrac{\partial\gamma}{\partial x}\left(  t,\zeta\right)
Z\left(  t\right)  +\tfrac{\partial\gamma}{\partial\mu}\left(  t,\zeta\right)
\eta\left(  t\right)  )\nu(d\zeta)dt]\nonumber\\
&  +\mathbb{E}[%
%TCIMACRO{\tint _{0}^{T}}%
%BeginExpansion
{\textstyle\int_{0}^{T}}
%EndExpansion
\langle p_{1}^{1}(t),DM^{\prime}(t)\rangle dt+%
%TCIMACRO{\tint _{0}^{T}}%
%BeginExpansion
{\textstyle\int_{0}^{T}}
%EndExpansion
DM(t)dp_{1}^{1}(t)]\nonumber\\
&  =\mathbb{E}[%
%TCIMACRO{\tint _{0}^{T}}%
%BeginExpansion
{\textstyle\int_{0}^{T}}
%EndExpansion
p_{1}^{0}(t)(\tfrac{\partial b}{\partial x}\left(  t\right)  Z\left(
t\right)  +\tfrac{\partial b}{\partial\mu}\left(  t\right)  \eta\left(
t\right)  )dt-%
%TCIMACRO{\tint _{0}^{T}}%
%BeginExpansion
{\textstyle\int_{0}^{T}}
%EndExpansion
\tfrac{\partial H_{1}}{\partial x}(t)Z(t)dt\nonumber\\
&  +%
%TCIMACRO{\tint _{0}^{T}}%
%BeginExpansion
{\textstyle\int_{0}^{T}}
%EndExpansion
q_{1}^{0}(t)(\tfrac{\partial\sigma}{\partial x}\left(  t\right)  Z\left(
t\right)  +\tfrac{\partial\sigma}{\partial\mu}\left(  t\right)  \eta\left(
t\right)  )dt\nonumber\\
&  +%
%TCIMACRO{\tint _{0}^{T}}%
%BeginExpansion
{\textstyle\int_{0}^{T}}
%EndExpansion%
%TCIMACRO{\tint _{\mathbb{R}_{0}}}%
%BeginExpansion
{\textstyle\int_{\mathbb{R}_{0}}}
%EndExpansion
r_{1}^{0}(t,\zeta)(\tfrac{\partial\gamma}{\partial x}\left(  t,\zeta\right)
Z\left(  t\right)  +\tfrac{\partial\gamma}{\partial\mu}\left(  t,\zeta\right)
\eta\left(  t\right)  )\nu(d\zeta)dt\nonumber\\
&  +{\textstyle\int_{0}^{T}}\langle p_{1}^{1}(t),DM^{\prime}(t)\rangle
dt-{\textstyle\int_{0}^{T}}\langle\nabla_{m}H_{1}(t),DM(t)\rangle dt].
\end{align}
Combining the above and recalling that $\eta$ is of the form (\ref{eta}), we
conclude that%

\[
0=\mathbb{E}[%
%TCIMACRO{\tint _{0}^{T}}%
%BeginExpansion
{\textstyle\int_{0}^{T}}
%EndExpansion
\tfrac{\partial H_{1}}{\partial\mu}(t)\eta(t)dt]=\mathbb{E}[%
%TCIMACRO{\tint _{s}^{T}}%
%BeginExpansion
{\textstyle\int_{s}^{T}}
%EndExpansion
\tfrac{\partial H_{1}}{\partial\mu}(t)\alpha_{1}dt]\text{; }s\geq t_{0}.
\]
Differentiating with respect to $s$ we obtain%
\begin{align*}
0  &  =\mathbb{E}[\tfrac{\partial H_{1}}{\partial\mu}(s)\alpha_{1}]\\
&  =\mathbb{E}[\tfrac{\partial H_{1}}{\partial\mu}(t_{0})|\mathcal{G}_{t_{0}%
}^{(1)}],
\end{align*}
because this holds for all $\alpha_{1}$ and all $s\geq t_{0}$.\newline This
argument can be reversed, to prove that (ii)$\Longrightarrow$(i). We omit the
details.\newline In the same manner, we can get the equivalence between
\[
\tfrac{d}{ds}J_{2}(\mu,u+s\pi)|_{s=0}=0
\]
and
\[
\mathbb{E}[\tfrac{\partial H_{2}}{\partial u}(t)|\mathcal{G}_{t}^{(2)}]=0.
\]
\hfill$\square$ \newline In the next section we will consider the zero-sum
case, and find conditions for a saddle point of such games.

\subsection{Zero-sum game}

\noindent In this section, we proceed to study the maximum principle for the
zero-sum game case. Let us then define the\ performance functional as
\[
J\left(  \mu,u\right)  =\mathbb{E[}g\left(  X(T),M(T)\right)  +%
%TCIMACRO{\tint _{0}^{T}}%
%BeginExpansion
{\textstyle\int_{0}^{T}}
%EndExpansion
\ell\left(  s,X\left(  s\right)  ,M(s),\mu(s),u\left(  s\right)  \right)
ds],
\]
where the state $X(t)$ is the solution of a SDE $\left(  \ref{sde}\right)
.$\newline The functions
\[
\ell(s,x,m,\mu,u)=\ell(s,x,m,\mu,u,\omega):\left[  0,T\right]  \times
\mathbb{R}\times\mathcal{M}_{0}\times\mathcal{M}\times\mathcal{U}\times
\Omega\rightarrow\mathbb{R}%
\]
and
\[
g(x,m)=g(x,m,\omega):\mathbb{R}\times\mathcal{M}_{0}\times\Omega
\rightarrow\mathbb{R}%
\]
are supposed to satisfy the following conditions:\newline

\begin{description}
\item[(a)] $\ell$ and $g$ are continuously differentiable with respect to
$x,u$ and admits Fréchet derivatives with respect to $m$ and $\mu$.

\item[(b)] Moreover, the function
\[%
%TCIMACRO{\U{211d} }%
%BeginExpansion
\mathbb{R}
%EndExpansion
\times\mathcal{M}_{0}\ni(x,m)\mapsto g(x,m)
\]

\end{description}

\noindent is required to be affine $\mathbb{P}$-a.s. \newline We consider the
stochastic zero-sum game to find $(\mu^{\ast},u^{\ast})$ such that%

\[
\underset{u\in\mathcal{A}_{\mathbb{G}}}{\sup}\underset{\mu\in\mathbb{M}%
_{\mathbb{G}}}{\inf}J(\mu,u)=\underset{\mu\in\mathbb{M}_{\mathbb{G}}}{\inf
}\underset{u\in\mathcal{A}_{\mathbb{G}}}{\sup}J(\mu,u)=J(\mu^{\ast},u^{\ast
}).
\]
We call $(\mu^{\ast},u^{\ast})$ a \textit{saddle point} for $J(\mu
,u).$\newline In this case, let the Hamiltonian
\[
H:[0,T]\times\mathbb{R}\times\mathcal{M}_{0} \times\mathcal{M}\times
\mathcal{U}\times\mathbb{R}\times\mathbb{R}\times\mathcal{R}\times
C_{a}([0,T],\mathcal{M}_{0})\rightarrow\mathbb{R}%
\]
be given by%

\begin{align*}
H(t,x,m, \mu, p^{0},q^{0},r^{0}(\cdot),p^{1})  &  =\ell(t,x,m,\mu
,u)+p^{0}b(t,x,\mu,u)+q^{0}\sigma(t,x,\mu,u)\\
&  +%
%TCIMACRO{\tint _{\mathbb{R}_{0}}}%
%BeginExpansion
{\textstyle\int_{\mathbb{R}_{0}}}
%EndExpansion
r^{0}(\zeta)\gamma\left(  t,x,\mu,u,\zeta\right)  \nu(d\zeta)+\langle
p^{1},\beta(m)\rangle.
\end{align*}
We assume the following:

\begin{description}
\item[(c)] $H$ is continuously differentiable with respect to $x,u$ and admits
Fréchet derivatives with respect to $m$ and $\mu$.

\item[(d)] The Hamiltonian function
\[
\mathbb{R}\times\mathcal{M}_{0}\times\mathcal{M}\times\mathcal{U}\ni
(x,m,\mu,u)\mapsto H(t,x,m,\mu, p^{0},q^{0},r^{0}(\cdot),p^{1})
\]
is \emph{convex} with respect to $(x,m,\mu)$ and \emph{concave} with respect
to $(x,m,u)$ $\mathbb{P}$.a.s and for each $t\in\left[  0,T\right]  $ ,
$p^{0},q^{0},r^{0}(\cdot)$ and $p^{1}.$\newline For $u\in\mathcal{A}%
_{\mathbb{G}},\mu\in\mathbb{M}_{\mathbb{G}}$ with corresponding solution
$X=X^{\mu,u}$, define $p=p^{\mu,u},q=q^{\mu,u}$ and $r=r^{\mu,u\text{ }}$ by
the adjoint equations: the real-BSDE in the unknown $(p^{0},q^{0},r^{0}%
)\in\mathcal{S}^{2}\times\mathbb{L}^{2}\times\mathbb{L}_{\nu}^{2}$ has the
following form
\end{description}

\begin{equation}
\left\{
\begin{array}
[c]{ll}%
dp^{0}(t) & =-\tfrac{\partial H}{\partial x}\left(  t\right)  dt+q^{0}%
(t)dB(t)+%
%TCIMACRO{\tint _{\mathbb{R}_{0}}}%
%BeginExpansion
{\textstyle\int_{\mathbb{R}_{0}}}
%EndExpansion
r^{0}(t,\zeta)\tilde{N}(dt,d\zeta)\text{; }t\in\left[  0,T\right]  \text{,}\\
p^{0}(T) & =\tfrac{\partial g}{\partial x}(X(T),M(T)),
\end{array}
\right.  \label{pro}%
\end{equation}

\noindent and the operator-valued BSDE for the unknown $(p^{1},q^{1},r^{1}%
)\in\mathcal{S}_{\mathbb{K}}^{2}\times\mathbb{L}_{\mathbb{K}}^{2}%
\times\mathbb{L}_{\nu,\mathbb{K}}^{2}$ is given by
\begin{equation}
\left\{
\begin{array}
[c]{lll}%
dp^{1}(t) & = & -\nabla_{m}H(t)dt+q^{1}(t)dB(t)+%
%TCIMACRO{\tint _{\mathbb{R}_{0}}}%
%BeginExpansion
{\textstyle\int_{\mathbb{R}_{0}}}
%EndExpansion
r^{1}(t,\zeta)\tilde{N}(dt,d\zeta);\quad t\in\lbrack0,T],\\
p^{1}(T) & = & \nabla_{m}g(X(T),M(T)).
\end{array}
\right.  \label{pro1}%
\end{equation}

\begin{theorem}
[Sufficient zero-sum maximum principle]Let $\left(  \hat{\mu},\hat{u}\right)
\in\mathbb{M}_{\mathbb{G}}\mathcal{\times A}_{\mathbb{G}}$ with corresponding
solutions $\hat{X}$ and $(p^{0},q^{0},r^{0})$, $(p^{1},q^{1},r^{1})$ of the
forward and backward stochastic differential equations $\left(  \ref{sde}%
\right)  ,\left(  \ref{pro}\right)  -\left(  \ref{pro1}\right)  ,$
respectively. Assume the following:

\begin{itemize}
\item
\[
\mathbb{E}[\hat{H}(t)|\mathcal{G}_{t}^{(1)}]=\underset{\mu\in\mathbb{M}%
_{\mathbb{G}}}{ess\text{ }\sup}\mathbb{E}[\check{H}(t)|\mathcal{G}_{t}%
^{(1)}],
\]

\item
\[
\mathbb{E}[\bar{H}(t)|\mathcal{G}_{t}^{(2)}]=\underset{u\in\mathcal{A}%
_{\mathbb{G}}}{ess\text{ }\sup}\mathbb{E}[\breve{H}(t)|\mathcal{G}_{t}%
^{(2)}],
\]

\end{itemize}

$\mathbb{P}$- a.s and for all $t\in\left[  0,T\right]  ,$ and that assumptions
(a)-(d) hold.\newline

Then $\left(  \hat{\mu},\hat{u}\right)  $ is a saddle point for $J\left(
\mu,u\right)  $.
\end{theorem}

\noindent This result will be applied in the next section.

\begin{theorem}
[Necessary zero-sum maximum principle]Let $\left(  \hat{\mu},\hat{u}\right)
\in\mathbb{M}\mathcal{_{\mathbb{G}}\times A}_{\mathbb{G}}$ with corresponding
solutions $\hat{X}$, $(p_{i}^{0},q_{i}^{0},r_{i}^{0})$ and $(p_{i}^{1}%
,q_{i}^{1},r_{i}^{1})$ of the forward and the backward stochastic differential
equations $\left(  \ref{sde}\right)  $ and $\left(  \ref{pro}\right)  -\left(
\ref{pro1}\right)  $, respectively, with corresponding derivative process
$\hat{Z}$ given by $\left(  \ref{dervz}\right)  .$ Then we have equivalence
between%
\[
\tfrac{d}{d\lambda}J(\mu+\lambda\eta,u)|_{\lambda=0}=\tfrac{d}{ds}J(\mu
,u+s\pi)|_{s=0}=0,
\]
and%
\[
\mathbb{E}[\tfrac{\partial H}{\partial\mu}(t)|\mathcal{G}_{t}^{(1)}%
]=\mathbb{E}[\tfrac{\partial H}{\partial u}(t)|\mathcal{G}_{t}^{(2)}]=0.
\]

\end{theorem}

\noindent{Proof.} \quad The same proof of both the sufficient and the
necessary maximum principles for the nonzero-sum games works for the zero-sum
case.\hfill$\square$ \newline

\section{Optimal consumption of a mean-field cash flow under uncertainty}

\noindent Consider a net cash flow $X^{\mu,\rho}=X$ modeled by%

\[
\left\{
\begin{array}
[c]{l}%
dX(t)=\left[  \mu(t)(V)-\rho(t)\right]  X(t)dt+\sigma\left(  t\right)
X(t)dB(t)+\int_{%
%TCIMACRO{\U{211d} }%
%BeginExpansion
\mathbb{R}
%EndExpansion
_{0}}\gamma\left(  t,\zeta\right)  X(t)\tilde{N}(dt,d\zeta)\text{; }%
t\in\left[  0,T\right]  \text{,}\\
X\left(  0\right)  =x>0\text{,}%
\end{array}
\right.
\]
where $\rho(t)\geq0$ is our \emph{relative consumption rate} at time $t$,
assumed to be a càdlàg, $\mathcal{G}_{t}^{(2)}$-adapted process. Here $V$ is a
given Borel subset of $\mathbb{R}$. The value of $\mu(t)$ on $V$ models the
relative growth rate of the cash flow. The relative consumption rate $\rho(t)$
is our control process. We assume that $%
%TCIMACRO{\tint _{0}^{T}}%
%BeginExpansion
{\textstyle\int_{0}^{T}}
%EndExpansion
\rho(t)dt<\infty$ a.s. This implies that $X(t)>0$ for all $t$, a.s. However,
the measure-valued process $\mu(t)$ represents a kind of scenario uncertainty,
and we want to maximise the total expected utility of the relative consumption
rate $\rho$ in the worst possible scenario $\mu$. We penalize $\mu(\cdot)$ for
being far away from the law process $\mathcal{L}(X(\cdot))$, in the sense that
we introduce a quadratic cost rate $[(\mu(t)-M(t))(V)]^{2}$ in the performance
functional. Hence we consider the zero-sum game%

\[
\underset{\rho}{\sup\text{ }}\underset{\mu}{\inf}\text{ }\mathbb{E}[%
%TCIMACRO{\tint _{0}^{T}}%
%BeginExpansion
{\textstyle\int_{0}^{T}}
%EndExpansion
\{\log(\rho(t)X(t))+[(\mu(t)-M(t))(V)]^{2}\}dt+\theta\log(X(T))],
\]
where $\theta=\theta(\omega)>0$ is a given bounded $\mathcal{F}_{T}%
$-measurable random variable, expressing the importance of the terminal value
$X(T)$. \textit{Here we have chosen a logarithmic utility because it is a
central choice, and in many cases, as here, this leads to a nice explicit
solution of the corresponding control problem.}\newline The Hamiltonian for
this zero-sum game takes the form%
\begin{align*}
H(t)  &  =\log(\rho x)+(\mu(V)-m(V))^{2}+p^{0}[\mu(V)x-\rho x]+q^{0}%
\sigma(t)x\\
&  \text{ \ \ \ \ \ \ \ \ \ \ \ \ }+%
%TCIMACRO{\tint _{\mathbb{R}_{0}}}%
%BeginExpansion
{\textstyle\int_{\mathbb{R}_{0}}}
%EndExpansion
r^{0}(\zeta)\gamma(t,\zeta)x\nu(d\zeta)+\langle p^{1},\beta(m)\rangle,
\end{align*}
and the adjoint processes $(p^{0},q^{0},r^{0})\in\mathcal{S}^{2}%
\times\mathbb{L}^{2}\times\mathbb{L}_{\nu}^{2}$,$(p^{1},q^{1},r^{1}%
)\in\mathcal{S}_{\mathbb{K}}^{2}\times\mathbb{L}_{\mathbb{K}}^{2}%
\times\mathbb{L}_{\nu,\mathbb{K}}^{2}$ are given by the BSDEs

\begin{itemize}
\item
\[
\left\{
\begin{array}
[c]{lll}%
dp^{0}(t) & = & -[\tfrac{1}{X(t)}+p^{0}(t)[\mu(t)(V)-\rho(t)]+q^{0}%
(t)\sigma(t)+%
%TCIMACRO{\tint _{\mathbb{R}_{0}}}%
%BeginExpansion
{\textstyle\int_{\mathbb{R}_{0}}}
%EndExpansion
r^{0}(t,\zeta)\gamma(t,\zeta)\nu(d\zeta)]dt\\
&  & +q^{0}(t)dB(t)+%
%TCIMACRO{\tint _{\mathbb{R}_{0}}}%
%BeginExpansion
{\textstyle\int_{\mathbb{R}_{0}}}
%EndExpansion
r^{0}(t,\zeta)\tilde{N}(dt,d\zeta)\text{;}\quad t\in\left[  0,T\right]  ,\\
p^{0}(T) & = & \tfrac{\theta}{X(T)}\text{,}%
\end{array}
\right.
\]

\item
\[
\left\{
\begin{array}
[c]{ll}%
dp^{1}(t) & =-\{2[\hat{\mu}(t)(V)-\hat{M}(t)(V)]\chi_{V}(\cdot)+<p^{1}%
(t),\beta(\cdot)>\}dt+q^{1}(t)dB(t)\\
& \text{ \ \ }+\int_{%
%TCIMACRO{\U{211d} }%
%BeginExpansion
\mathbb{R}
%EndExpansion
_{0}}r^{1}(t,\zeta)\tilde{N}(dt,d\zeta);\quad t\in\lbrack0,T],\\
p^{1}(T) & =0,
\end{array}
\right.
\]

\end{itemize}

\noindent where $\chi_{V}(\cdot)$ is the operator which evaluates a given
measure at $V$, i.e. $\left\langle \chi_{V},\lambda\right\rangle =\lambda(V)$
for all $\lambda\in\mathcal{M}_{0}$. The first order condition for the optimal
consumption rate $\hat{\rho}$ is%

\[
\mathbb{E}[\tfrac{1}{\hat{\rho}(t)}-\hat{p}^{0}(t)\hat{X}(t)|\mathcal{G}%
_{t}^{(2)}]=0.
\]
Since $\hat{\rho}(t)$ is $\mathcal{G}_{t}^{(2)}$-adapted, we have
\[
\hat{\rho}(t)=\tfrac{1}{\mathbb{E}[\hat{p}^{0}(t)\hat{X}(t)|\mathcal{G}%
_{t}^{(2)}]}.
\]
Now we use the minimum condition with respect to $\mu$ at $\mu=\hat{\mu}$ and
get
\[
\mathbb{E}[2[\hat{\mu}(t)(V)-\hat{M}(t)(V)]\lambda(V)+\hat{p}^{0}(t)\hat
{X}(t)\lambda(V)|\mathcal{G}_{t}^{(1)}]=0\text{, for all }\lambda
\in\mathcal{M}_{0}.
\]
Using that $\hat{\mu}(t)$ is $\mathcal{G}_{t}^{(1)}$-adapted, we obtain
\[
\hat{\mu}(t)(V)=\mathbb{E}[\hat{M}(t)(V)-\tfrac{1}{2}\hat{p}^{0}(t)\hat
{X}(t)|\mathcal{G}_{t}^{(1)}].
\]
It remains to find $\hat{p}^{0}(t)\hat{X}(t)$: We have by applying the Itô
formula to $P(t):=\hat{p}^{0}(t)\hat{X}(t)$:%

\begin{align}
dP(t)  &  =\hat{p}^{0}(t)d\hat{X}(t)+\hat{X}(t)d\hat{p}^{0}(t)+d[\hat{p}%
^{0},\hat{X}]_{t}\nonumber\\
&  =\hat{p}^{0}(t)([\left(  \hat{\mu}(t)(V)-\rho(t)\right)  \hat{X}%
(t)]dt+\hat{\sigma}\left(  t\right)  \hat{X}(t)dB(t)+%
%TCIMACRO{\tint _{\mathbb{R}_{0}}}%
%BeginExpansion
{\textstyle\int_{\mathbb{R}_{0}}}
%EndExpansion
\hat{\gamma}\left(  t,\zeta\right)  \hat{X}(t)\tilde{N}(dt,d\zeta))\nonumber\\
&  +\hat{X}(t)[-\tfrac{1}{\hat{X}(t)}-\hat{p}^{0}(t)[\hat{\mu}(t)(V)-\rho
(t)]-\hat{q}^{(0)}(t)\sigma(t)-%
%TCIMACRO{\tint _{\mathbb{R}_{0}}}%
%BeginExpansion
{\textstyle\int_{\mathbb{R}_{0}}}
%EndExpansion
\hat{r}^{0}(t,\zeta)\hat{\gamma}(t,\zeta)\nu(d\zeta)]dt\nonumber\\
&  +\hat{q}^{0}(t)\hat{X}(t)dB(t)+%
%TCIMACRO{\tint _{\mathbb{R}_{0}}}%
%BeginExpansion
{\textstyle\int_{\mathbb{R}_{0}}}
%EndExpansion
\hat{r}^{0}(t,\zeta)\hat{X}(t)\tilde{N}(dt,d\zeta)+\hat{q}^{0}(t)\hat{\sigma
}\left(  t\right)  \hat{X}(t)dt\nonumber\\
&  +%
%TCIMACRO{\tint _{\mathbb{R}_{0}}}%
%BeginExpansion
{\textstyle\int_{\mathbb{R}_{0}}}
%EndExpansion
\hat{r}^{0}(t,\zeta)\hat{\gamma}(t,\zeta)\hat{X}(t)N(dt,d\zeta). \label{na1}%
\end{align}
By definition
\begin{equation}%
\begin{array}
[c]{ll}%
%TCIMACRO{\tint _{\mathbb{R}_{0}}}%
%BeginExpansion
{\textstyle\int_{\mathbb{R}_{0}}}
%EndExpansion
\hat{r}^{0}(t,\zeta)\hat{\gamma}(t,\zeta)\hat{X}(t)\tilde{N}(dt,d\zeta) & =%
%TCIMACRO{\tint _{\mathbb{R}_{0}}}%
%BeginExpansion
{\textstyle\int_{\mathbb{R}_{0}}}
%EndExpansion
\hat{r}^{0}(t,\zeta)\hat{\gamma}(t,\zeta)\hat{X}(t)N(dt,d\zeta)\\
& \text{ \ }-%
%TCIMACRO{\tint _{\mathbb{R}_{0}}}%
%BeginExpansion
{\textstyle\int_{\mathbb{R}_{0}}}
%EndExpansion
\hat{r}^{0}(t,\zeta)\hat{\gamma}(t,\zeta)\hat{X}(t)\nu(d\zeta)dt.
\end{array}
\label{na2}%
\end{equation}
Substituting $\left(  \ref{na2}\right)  $\ in $\left(  \ref{na1}\right)  $ yields%

\[%
\begin{array}
[c]{ll}%
dP(t) & =-dt+[P(t)\hat{\sigma}(t)+\hat{q}^{0}(t)\hat{X}(t)]dB(t)\\
& +%
%TCIMACRO{\tint _{\mathbb{R}_{0}}}%
%BeginExpansion
{\textstyle\int_{\mathbb{R}_{0}}}
%EndExpansion
[P(t)\hat{\gamma}(t,\zeta)+\hat{r}^{0}(t,\zeta)\hat{X}(t)(1+\hat{\gamma
}(t,\zeta))]\tilde{N}(dt,d\zeta).
\end{array}
\]
Hence, if we put
\[%
\begin{array}
[c]{lll}%
P(t) & := & \hat{p}^{0}(t)\hat{X}(t),\\
Q(t) & := & P(t)\hat{\sigma}(t)+\hat{X}(t)\hat{q}^{0}(t),\\
R(t,\zeta) & := & P(t)\hat{\gamma}(t,\zeta)+\hat{r}^{0}(t,\zeta)\hat
{X}(t)(1+\hat{\gamma}(t,\zeta)).
\end{array}
\]
with $(P,Q,R)\in\mathcal{S}^{2}\times\mathbb{L}^{2}\times\mathbb{L}_{\nu}^{2}$
satisfies the BSDE%

\[
\left\{
\begin{array}
[c]{lll}%
dP(t) & = & -dt+Q(t)dB(t)+%
%TCIMACRO{\tint _{\mathbb{R}_{0}}}%
%BeginExpansion
{\textstyle\int_{\mathbb{R}_{0}}}
%EndExpansion
R(t,\zeta)\tilde{N}(dt,d\zeta);\text{ \ \ }t\in\left[  0,T\right]  ,\\
P(T) & = & \theta.
\end{array}
\right.
\]
Solving this BSDE as in $\left(  \ref{solutio}\right)  $, we find the closed
formula for $P(t)$ as
\[%
\begin{array}
[c]{lll}%
P(t) & = & \mathbb{E}[\theta+\int_{t}^{T}ds|\mathcal{F}_{t}]\\
& = & \mathbb{E}\left[  \theta|\mathcal{F}_{t}\right]  +T-t.
\end{array}
\]
Hence we have proved the following:

\begin{theorem}
The optimal consumption rate $\hat{\rho}(t)$ and the optimal model uncertainty
law $\hat{\mu}(t)$ are given respectively in feed-back form by
\[%
\begin{array}
[c]{lll}%
\hat{\rho}(t) & = & \tfrac{1}{T-t+\mathbb{E}[\theta|\mathcal{G}_{t}^{(2)}]},\\
\hat{\mu}(t)(V) & = & \hat{M}(t)(V)+T-t-\tfrac{1}{2}\mathbb{E}[\theta
|\mathcal{G}_{t}^{(1)}].
\end{array}
\]

\end{theorem}

\section{Appendix}

\noindent Let us give now the rest of the proof of Theorem \ref{Thm}. We want
to prove that $J_{2}(\hat{\mu},u)\leq J_{2}(\hat{\mu},\hat{u})$. Using
definition $\left(  \ref{perf}\right)  $ gives for fixed $\hat{\mu}%
\in\mathbb{M}_{\mathbb{G}}$ and an arbitrary $u\in\mathcal{A}_{\mathbb{G}}$
\begin{equation}%
\begin{array}
[c]{lll}%
J_{2}(\hat{\mu},u)-J_{2}(\hat{\mu},\hat{u}) & = & j_{1}+j_{2},
\end{array}
\label{J2}%
\end{equation}
where
\begin{align*}
&
\begin{array}
[c]{lll}%
j_{1} & = & \mathbb{E}[\int_{0}^{T}\left\{  \breve{\ell}_{2}(t)-\bar{\ell}%
_{2}(t)\right\}  dt],
\end{array}
\\
&
\begin{array}
[c]{lll}%
j_{2} & = & \mathbb{E}[\breve{g}_{2}(X(T),M(T))-\bar{g}_{2}(\hat{X}(T),\hat
{M}(T))].
\end{array}
\end{align*}
Applying the definition of the Hamiltonian $\left(  \ref{haml}\right)  $ we
have%
\begin{equation}%
\begin{array}
[c]{l}%
j_{1}=\mathbb{E}[%
%TCIMACRO{\tint _{0}^{T}}%
%BeginExpansion
{\textstyle\int_{0}^{T}}
%EndExpansion
\{\breve{H}_{2}(t)-\breve{H}_{2}(t)-\hat{p}_{2}^{0}(t)\tilde{b}(t)-\hat{q}%
_{2}^{0}(t)\tilde{\sigma}(t)\\
-{%
%TCIMACRO{\tint _{\mathbb{R}_{0}}}%
%BeginExpansion
{\textstyle\int_{\mathbb{R}_{0}}}
%EndExpansion
}\hat{r}_{2}^{0}(t,\zeta)\tilde{\gamma}(t,\zeta)\nu(d\zeta)-\langle\hat{p}%
_{2}^{1}(t),\tilde{M}^{\prime}(t)\rangle\}dt],
\end{array}
\label{j1}%
\end{equation}
where $\tilde{b}(t)=\breve{b}(t)-\bar{b}(t)$. etc., and
\[
\tilde{M}^{\prime}(t)=\tfrac{d\tilde{M}(t)}{dt}.
\]
Concavity of $g_{2}$ and the definition of the terminal value of the BSDEs
$\left(  \ref{eqp0}\right)  $ and $\left(  \ref{eqp1}\right)  $ shows that%
\begin{align}
j_{2}  &  \leq\mathbb{E}[\tfrac{\partial g_{2}}{\partial x}(T)\tilde
{X}(T)+\langle\nabla_{m}g_{2}(T),\tilde{M}(t)\rangle]\nonumber\\
&  =\mathbb{E}[\hat{p}_{2}^{0}(T)\tilde{X}(T)+\langle\hat{p}_{2}^{1}%
(T),\tilde{M}(t)\rangle]. \label{j2}%
\end{align}
Applying the Itô formula to $\hat{p}_{2}^{0}\tilde{X}$ and $\langle\hat{p}%
_{2}^{1},\tilde{M}\rangle$, we get
\begin{align*}
j_{2}  &  \leq\mathbb{E}[\hat{p}_{2}^{0}(T)\tilde{X}(T)+\langle\hat{p}_{2}%
^{1}(T),\tilde{M}(T)\rangle]\\
&  =\mathbb{E}[%
%TCIMACRO{\tint _{0}^{T}}%
%BeginExpansion
{\textstyle\int_{0}^{T}}
%EndExpansion
\hat{p}_{2}^{0}(t)d\tilde{X}(t)+%
%TCIMACRO{\tint _{0}^{T}}%
%BeginExpansion
{\textstyle\int_{0}^{T}}
%EndExpansion
\tilde{X}(t)d\hat{p}_{2}^{0}(t)+%
%TCIMACRO{\tint _{0}^{T}}%
%BeginExpansion
{\textstyle\int_{0}^{T}}
%EndExpansion
\hat{q}_{2}^{0}(t)\tilde{\sigma}(t)dt+%
%TCIMACRO{\tint _{0}^{T}}%
%BeginExpansion
{\textstyle\int_{0}^{T}}
%EndExpansion%
%TCIMACRO{\tint _{\mathbb{R}_{0}}}%
%BeginExpansion
{\textstyle\int_{\mathbb{R}_{0}}}
%EndExpansion
\hat{r}_{2}^{0}(t,\zeta)\tilde{\gamma}(t,\zeta)\nu(d\zeta)dt]\\
&  +\mathbb{E}[%
%TCIMACRO{\tint _{0}^{T}}%
%BeginExpansion
{\textstyle\int_{0}^{T}}
%EndExpansion
\langle\hat{p}_{2}^{1}(t),d\tilde{M}(t)\rangle+%
%TCIMACRO{\tint _{0}^{T}}%
%BeginExpansion
{\textstyle\int_{0}^{T}}
%EndExpansion
\tilde{M}(t)d\tilde{p}_{2}^{1}(t)]\\
&  =\mathbb{E}[%
%TCIMACRO{\tint _{0}^{T}}%
%BeginExpansion
{\textstyle\int_{0}^{T}}
%EndExpansion
\hat{p}_{2}^{0}(t)\tilde{b}(t)dt-%
%TCIMACRO{\tint _{0}^{T}}%
%BeginExpansion
{\textstyle\int_{0}^{T}}
%EndExpansion
\tfrac{\partial\bar{H}_{2}}{\partial x}(t)\tilde{X}(t)dt+%
%TCIMACRO{\tint _{0}^{T}}%
%BeginExpansion
{\textstyle\int_{0}^{T}}
%EndExpansion
\hat{q}_{2}^{0}(t)\tilde{\sigma}(t)dt\\
&  +%
%TCIMACRO{\tint _{0}^{T}}%
%BeginExpansion
{\textstyle\int_{0}^{T}}
%EndExpansion%
%TCIMACRO{\tint _{\mathbb{R}_{0}}}%
%BeginExpansion
{\textstyle\int_{\mathbb{R}_{0}}}
%EndExpansion
\hat{r}_{2}^{0}(t,\zeta)\tilde{\gamma}(t,\zeta)\nu(d\zeta)dt+%
%TCIMACRO{\tint _{0}^{T}}%
%BeginExpansion
{\textstyle\int_{0}^{T}}
%EndExpansion
\langle\hat{p}_{2}^{1}(t),\tilde{M}^{\prime}(t)\rangle dt-%
%TCIMACRO{\tint _{0}^{T}}%
%BeginExpansion
{\textstyle\int_{0}^{T}}
%EndExpansion
\langle\nabla_{m}\bar{H}_{2}(t),\tilde{M}(t)\rangle dt],
\end{align*}
where we have used that the $dB(t)$ and $\tilde{N}(dt,d\zeta)$ integrals have
mean zero. Substituting $\left(  \ref{j1}\right)  $ and $\left(
\ref{j2}\right)  $ into $\left(  \ref{J2}\right)  $, we obtain%
\[
J_{2}(\hat{\mu},u)-J_{2}(\hat{\mu},\hat{u})\leq\mathbb{E}[%
%TCIMACRO{\tint _{0}^{T}}%
%BeginExpansion
{\textstyle\int_{0}^{T}}
%EndExpansion
\{\breve{H}_{2}(t)-\bar{H}_{2}(t)-\tfrac{\partial\bar{H}_{2}}{\partial
x}(t)\tilde{X}(t)-\langle\nabla_{m}\bar{H}_{2}(t),\tilde{M}(t)\rangle\}dt].
\]
Since $H_{2}$ is concave and the process $u$ is $\mathcal{G}_{t}^{(2)}%
$-adapted, we have%
\begin{align}
J_{2}(\hat{\mu},u)-J_{2}(\hat{\mu},\hat{u})  &  \leq\mathbb{E}[%
%TCIMACRO{\tint _{0}^{T}}%
%BeginExpansion
{\textstyle\int_{0}^{T}}
%EndExpansion
\tfrac{\partial\bar{H}_{2}}{\partial u}(t)\left(  u(t)-\hat{u}(t)\right)
dt]\nonumber\\
&  =\mathbb{E}[%
%TCIMACRO{\tint _{0}^{T}}%
%BeginExpansion
{\textstyle\int_{0}^{T}}
%EndExpansion
\mathbb{E[}\tfrac{\partial\bar{H}_{2}}{\partial u}(t)|\mathcal{G}_{t}%
^{(2)}]\left(  u(t)-\hat{u}(t)\right)  dt]\nonumber\\
&  \leq0,\nonumber
\end{align}
because $\bar{H}_{2}$ has a maximum at $\hat{u}$.\hfill$\square$ \newline

\noindent$\mathbf{Acknowledgments}$\newline We are grateful to Boualem
Djehiche for helpful comments.

\end{document}